\documentclass[10pt,reqno]{amsart}
\usepackage{amsmath,amsfonts,amssymb,amscd,amsthm,amsbsy,bbm, epsf,calc,enumerate,color,graphicx,verbatim,cite}
\usepackage{color}
\usepackage{datetime}

\newcounter{hours}\newcounter{minutes}

\makeatletter
\newtheorem*{rep@theorem}{\rep@title}
\newcommand{\newreptheorem}[2]{%
\newenvironment{rep#1}[1]{%
 \def\rep@title{#2 \ref{##1}}%
 \begin{rep@theorem}}%
 {\end{rep@theorem}}}
\makeatother

\theoremstyle{theorem}
\newtheorem{definition}{Theorem}[section]
\newtheorem{thm}{Theorem}[section]
\newreptheorem{thm}{Theorem}

\newtheorem{lem}[thm]{Lemma}
\newreptheorem{lem}{Lemma}

\theoremstyle{definition}
\newtheorem{ex}[thm]{Example}

\newtheorem{DEF}[definition]{Definition}

\theoremstyle{remark}                  
\newtheorem{rem}[thm]{Remark}

\def\real{{\mathbb R}}
\def\rational{{\mathbb Q}}

\def\integer{{\mathbb Z}}

\def\e{\varepsilon}

\def\Tr{\textnormal{Tr}}

\DeclareMathOperator*{\osc}{osc}

\def\liminf{\mathop{\lim\,\inf}\limits}%
\def\limsup{\mathop{\lim\,\sup}\limits}%

\numberwithin{equation}{section}

\begin{document}
\title{Homogenization of the Oscillating Dirichlet Boundary Condition in General Domains}
\author{William M. Feldman }
\address{Dept. of Mathematics, UCLA.}
\email{wfeldman10@math.ucla.edu }

\begin{abstract}
We prove the homogenization of the Dirichlet problem for fully nonlinear elliptic operators with periodic oscillation in the operator and of the boundary condition for a general class of smooth bounded domains.  This extends the previous results of Barles and Mironescu \cite{BarlesMironescu12} in half spaces.  We show that homogenization holds despite a possible lack of continuity in the homogenized boundary data.  The proof is based on a comparison principle with partial Dirichlet boundary data which is of independent interest.

\end{abstract}

\maketitle

\section{Introduction}

We consider the homogenization of the following Dirichlet problem set in a bounded domain in $\Omega \subset \real^n$ with smooth boundary,
\begin{equation}\label{eqn: epsilon prob}
\left\{
\begin{array}{lll}
F_\varepsilon(D^2u^\varepsilon,x,\tfrac{x}{\varepsilon}) = 0 & \hbox{ for } & x \in \Omega \\ \\
u^\varepsilon(x) = g(x,\tfrac{x}{\varepsilon}) & \hbox{ for } & x \in \partial \Omega,
\end{array}\right.
\end{equation}
with H\"{o}lder continuous boundary data $g \in C^{\alpha}(\real^n\times \real^n)$ for some $\alpha \in (0,1)$.  The operator $(M,x,y) \to F_\varepsilon (M,x,y)$ is assumed to be uniformly elliptic, and both the operator and the boundary condition $g(x,y)$ are assumed to be $\integer^n$-periodic in the $y$ variable.  See Section \ref{subsec: Assumptions on the Operators} for the precise assumptions we make on the operators.  

\medskip

The homogenization of the problem \eqref{eqn: epsilon prob} when $g(x,y) = g(x)$ is independent of the fast variable is somewhat classical at this point and was done by Evans \cite{Evans92,Evans89}, the ideas originated in the unpublished paper of Lions, Papanicolaou and Varadhan \cite{LPV88} on the periodic homogenization of Hamilton-Jacobi equations.  There are too many works investigating the homogenization of the interior operators in various settings for us to list all of them here, we simply mention that more recently the interior homogenization has been shown in the very general setting of stationary ergodic random media by Caffarelli, Souganidis and Wang in \cite{CSW05}.  

\medskip

As far as oscillating boundary data is concerned there is less literature.  For divergence form operators the case of oscillating co-normal Neumann data is in the classical book of Bensoussan, Papanicolaou and Lions \cite{BPL78}.  For the case of non co-normal Neumann data and fully nonlinear operators there are a number of recent works.  The papers of Arisawa \cite{Arisawa03} and Barles, Da Lio, Lions and Souganidis \cite{BDLS08} identify the cell problem and show the homogenization in half-space type domains.  The general domain case for the Neumann problem was solved by Choi and Kim in \cite{ChoiKim12} (see also \cite{CKL12}).  They require that the operator $F_\e$ homogenizes to a rotationally invariant operator, in this situation the homogenized boundary data is continuous which guarantees the uniqueness of the limiting problem and thus stability of the $\e$-problem as $\e\to 0$.  This issue is explained in more detail in \cite{ChoiKim12}.  Our work is based on the ideas in \cite{ChoiKim12,CKL12}.  Indeed, one of the novelties of this article is to show the uniqueness and stability result for a general class of operators even though the homogenized boundary data may be discontinuous \cite{CharlieInwon}.  In particular we do not need to put any special assumptions on the operators besides the uniform ellipticity.  We believe that similar ideas should apply to the case of oscillatory Neumann data.

\medskip

 For the oscillating Dirichlet data case there are even less references.  For linear divergence form equations the homogenization in general domains was shown recently by Gerard-Varet and Masmoudi in \cite{GVM11,GVM12} where they are able to consider linear systems, see also \cite{LS12}.  The work of Barles and Mironescu \cite{BarlesMironescu12} deals with fully nonlinear non-divergence form elliptic operators when $\Omega$ is a half-space.  Barles and Mironescu identify the cell problem and the attendant difficulties in solving it under quite general assumptions on the operators.  We will be concerned with extending this result to general domains.  In order to clarify the main issues that arise from the general domain rather than from the cell problem we will avoid giving the most general assumptions on the operators.  We expect that our results should also hold in the generality at which \cite{BarlesMironescu12} show the cell problem can be solved, but this will be addressed in future work.

\medskip

 At least on the surface a difficulty of this problem, especially in contrast to the corresponding oscillating Neumann data problem, is a lack of equicontinuity for the $u^\e$.  The optimal estimate for the continuity of the $u^\e$ up to the boundary cannot be any better than the continuity of $g(x,\tfrac{x}{\e})$ which has oscillations of size $1$ at arbitrarily small scales in the limit $\e \to 0$.  In particular the homogenization of \eqref{eqn: epsilon prob}, if it occurs, must be happening outside of a `boundary layer' of width $o(1)$.  Keeping this in mind, let us call $\overline{F}$ to be the homogenized operator on the interior given by the results of \cite{Evans92,Evans89,CSW05}.  Our goal is to show that the solutions $u^\e$ of the $\e$-problem \eqref{eqn: epsilon prob} converge locally uniformly in $\Omega$ to $\overline{u}$ solving an equation of the form,
\begin{equation}\label{eqn: hom eqn 1}
\left\{
\begin{array}{lll}
\overline{F}(D^2\overline{u},x) = 0 & \hbox{ for } & x \in \Omega \\ \\
\overline{u}(x) = \overline{g}(x) & \hbox{ for } & x \in \partial \Omega.
\end{array}\right.
\end{equation}
As we will see, the homogenized boundary condition $\overline{g}$ will depend not only on the function $g$, but also on $\nu_x$, the inner normal to $\Omega$ at $x$ and on the $F_\e$.

 \medskip
 
 Let us describe heuristically the approach taken in \cite{BarlesMironescu12} to identifying the homogenized boundary condition.  The analysis of the problem \eqref{eqn: epsilon prob} proceeds by blowing up about points $x \in \partial \Omega$.  As in \cite{BarlesMironescu12}, the introduction of the localizations,
$$ v^\e(y) = u^\e(x_0+\e y), $$
at points $x_0 \in \partial \Omega$ in the limit as $\varepsilon$ is taken to $0$ leads us formally to the following ``cell problem" set in the half space $P(x_0,\nu):= \{y: (y-x_0)\cdot \nu >0\}$ where $\nu \in S^{n-1} = \{x \in \real^n: |x| =1\}$,
\begin{equation}\label{eqn: cell prob 1}
\left\{
\begin{array}{lll}
F(D^2v,y) = 0  & \hbox{ in } & P(x_0,\nu) \\ \\
v(y) = \psi(y) & \hbox{ on } & \partial P(x_0,\nu).
\end{array}\right.
\end{equation}
  Some appropriate assumptions on the $F_\e$ and $g$ will yield that $F$ is uniformly elliptic and that both $F$ and $\psi$ are $\integer^n$-periodic in $y$.  Then by a Liouville type result the limit $v(y+R\nu)$ as $R \to \infty$ will be a constant $\overline{\mu}(\psi,F,\nu)$, we will identify this as the homogenized boundary condition.  Then, at least formally, one expects that $|u^\e(x+\e R\nu) -v(R\nu)| \to 0$ as $\e \to 0$.  In particular, if we pull away from the boundary sufficiently far near $x_0$, we can recover the interior averaging.  This suggests that our homogenized boundary condition in \eqref{eqn: hom eqn 1} should be $\overline{g}(x) = \overline{\mu}(g(x,\cdot),F^x,\nu_x)$ where $\nu_x$ is the inner unit normal to $\Omega$ at $x$ and $F^x$ is an appropriate blow up limit of $F_\e$.
  
  \medskip
 
  The main difficulty that one has to deal with in extending to general domains is the \textit{rational boundary points} $x \in \partial \Omega$ such that the interior normal $\nu_x$ to $\Omega$ at $x$ is in $ \real \integer ^n$.   Exactly when $\nu \in \real \integer^n$ the cell problem \eqref{eqn: cell prob 1} does not yield a unique constant $\overline{\mu}$.  This corresponds to an a-priori lack of uniform control in $\nu$ of the rate of convergence of $v(y+R\nu) \to \overline{\mu}(\nu)$. Some stability of this rate of convergence in $\nu$ is necessary in order to show that the behavior of the cell problem governs the homogenization near irrational boundary points.  
  
  \medskip
  
    It turns out that $v(x+R\nu) \to \overline{\mu}(\nu)$ uniformly in $\nu \in S^{n-1}$ is equivalent to $\overline{\mu}(\psi,F,\cdot)$ having a unique continuous extension from the irrational directions $S^{n-1} \setminus \real \integer^n$ to all of $S^{n-1}$.  When $\overline{F}$ is a rotationally invariant operator one can show using the methods of \cite{CKL12,ChoiKim12} that such a continuous extension exists.  For general $F$ the homogenized boundary condition $\overline{\mu}$ is expected to have no such continuous extension \cite{CharlieInwon}, so we are left to consider weaker continuity conditions which still suffice to prove the homogenization.  
  
  \medskip
  
  One of the main contributions of this work is showing that the weak stability which one has `intrinsically' is sufficient for the homogenization to hold.  We find this stability by carefully analyzing the rates of convergence for the cell problem in terms of the discrepancy function, see Sections \ref{subsec: Rational and Irrational Directions} and \ref{sec: The Cell Problem}.  
  
  \medskip
  
  Let us call $\Gamma$ the set of rational boundary points,
  $$ \Gamma(\Omega) := \{ x \in \partial \Omega : \nu_x \in \real \integer^n\}.$$
  In light of the above discussion, the best continuity result we expect for general $F_\e$ is,
  \begin{equation}\label{eqn: cont at irrat}
   \overline{g} : \partial \Omega \to \real \ \hbox{ is continuous at } x \in \partial\Omega \setminus \Gamma,
   \end{equation}
and it turns out that $\partial \Omega \setminus \Gamma$ is also exactly where one is able to exploit the connection between the cell problem and the $\e$-problem to show that any subsequential locally uniform limit of the $u^\e$ satisfies the homogenized boundary data.  In order to show that all subsequential limits are the same the fundamental question becomes, when is there uniqueness for the problem:
\begin{equation}\label{eqn: hom eqn 2}
\left\{
\begin{array}{lll}
\overline{F}(D^2\overline{u},x) = 0 & \hbox{ for } & x \in \Omega \\ \\
\overline{u}(x) = \overline{g}(x) & \hbox{ for } & x \in \partial \Omega \setminus \Gamma(\Omega).
\end{array}\right.
\end{equation}
 We are not aware of any existing results in this direction, it is worth noting that this problem is related to the open question of whether there exists a boundary version of the Alexandroff-Bakelman-Pucci estimate.  We provide a resolution of this question when $\Gamma(\Omega)$ has sufficiently small Hausdorff dimension in the following Theorem, see Section \ref{subsec: Comparison with Partial Boundary Data} for a precise explication of all the assumptions,
\begin{thm}\label{thm: ABP type 1} \textup{[Theorem \ref{thm: ABP type}]}
Let $\overline{F}(M,y)$ be uniformly elliptic with constants $(\lambda,\Lambda)$ and satisfy appropriate continuity assumptions in $(M,y)$.  Suppose that $\Omega$ has the strict $\gamma$ exterior cone condition for some $0<\gamma \leq 1$.  There exists $\beta_0(n,\lambda,\Lambda,\gamma)>\max\{0,\tfrac{\lambda}{\Lambda}(n-1)-1\}$ such that if
$$\beta_0>\dim_{\mathcal{H}} \Gamma(\Omega) \geq 0,$$
then \eqref{eqn: hom eqn 2} has a unique bounded solution $\overline{u}$.
\end{thm}

\begin{rem}
The constraint on the Hausdorff dimension of $\Gamma(\Omega)$ allows us to construct a sequence of super solution barriers $\phi_j$ which are large on $\Gamma(\Omega)$, but converge to zero on compact subsets of $\Omega$.  Then for any bounded solution $v$ of \eqref{eqn: hom eqn 2} one can use the classical comparison principle to show $v \leq \overline{u}+\phi_j$, and letting $j \to \infty$ will yield the result.
\end{rem}

\medskip

Now we will give a statement of our main result. The results are not stated in their full generality in order that it be understandable at this point in the article.  For the exact statements one can refer to the main text, in particular the main result is in Theorem \ref{thm: dirichlet hom}.

\begin{thm}\label{thm: main thm}(Solution of the Cell Problem and Homogenization)
\begin{enumerate}[(i)]
\item \textup{$[$Lemmas \ref{lem: ID BC}, \ref{lem: mu properties} and \ref{lem: continuity of BC}$]$}  Let $v_{x_0,\nu}$ solve \eqref{eqn: cell prob 1} in $P(\nu,x_0)$, there exists a constant $\overline{\mu}(\psi,F,\nu)$ such that for all $\nu \in S^{n-1} \setminus \real \integer^n$,
$$ \sup_{x_0 \in \real^n} \sup_{y \in P(x_0,\nu)} |v_{x_0,\nu} (y+R\nu)-\overline{\mu}(\psi,F,\nu)| \to 0 \ \hbox{ as } \ R \to \infty.$$
Moreover $\overline{\mu}$ is continuous in all of its arguments at irrational directions $\nu \in S^{n-1} \setminus \real \integer^n$ (for the precise meaning of this refer to the referenced Lemmas).
\item \textup{$[$Theorem \ref{thm: dirichlet hom}$]$}  Let $u^\e$ be the solutions of the $\e$-problem \eqref{eqn: epsilon prob} and $\Omega$ a bounded domain with $C^2$ boundary satisfying the condition of Theorem \ref{thm: ABP type 1}.  Then the $u^\e$ converge locally uniformly in $\Omega$ to the unique solution $\overline{u}$ of \eqref{eqn: hom eqn 2} with,
$$ \overline{g}(x) = \overline{\mu}(g(x,\cdot),F^x,\nu_x). $$
\end{enumerate}
\end{thm}
\begin{rem}\label{rem: flat parts}
In particular the homogenization result holds whenever the set of rational boundary points $\Gamma$ is countable for general uniformly elliptic operators $F_\e$ satisfying the assumptions of Section \ref{subsec: Assumptions on the Operators}.  For example, by the inverse function theorem, this includes the classical assumptions of \cite{BPL78} for the homogenization of co-normal Neumann data, that $\partial\Omega$ has no flat parts in the sense that $D_\tau\nu_x$ (the gradient in the tangential directions at $x$) has rank $n-1$.
\end{rem}
\begin{rem}Here $F^x$ is a rescaling of the $F_\e$ at $x \in \partial \Omega$.  Under the assumptions we will put on the $F_\e$, the operators $F_\e$ themselves will not actually depend on $\e$ and one can simply take $F^x(M,y) = F_\e(M,x,y)$.  For a precise characterization of the $F^x$ see Lemma \ref{lem: scale one op}.
\end{rem}
\begin{rem}\label{rem: nonlinearity of mu}
It turns out that $\overline{\mu}(\psi,F,\nu)$ is just the average of $\psi$ over a unit cell when $F$ is linear and $\nu$ is an irrational direction, see Section \ref{subsec: The Nonlinearity of the Homogenized Boundary Condition} for the proof.  Since there is no dependence on the operator in the linear case one may wonder whether $\overline{\mu}$ is always just the average even when $F$ is nonlinear.  We show that in fact this is not the case in Section \ref{subsec: The Nonlinearity of the Homogenized Boundary Condition}, generically when $F$ is nonlinear $\overline{\mu}(\psi,F,\nu)$ depends on $F$ and is not linear in the $\psi$ argument.
\end{rem}
\begin{rem}\label{rem: cell prob}
We do not consider the case of divergence form operators here, for example,
\begin{equation*}
\left\{
\begin{array}{lll}
-\operatorname{div}(A(x,\e^{-1}x)Du^\e) = 0 & \hbox{ in } & \Omega \\
u^\e(x) = g(x,\e^{-1}x) & \hbox{ on } & \partial\Omega
\end{array}\right.
\end{equation*}
In this case the rescaled operator $F_x$ in the cell problem has the form,
$$F_x(D^2u,Du,y) = -\Tr(A(x,y)D^2u)-\operatorname{div}_yA(x,y) \cdot Du.$$
We believe that the only major difference to handle this kind of operator lies in the Localization Lemma \ref{lem: localization}.  As in \cite{BarlesMironescu12} the requirement is that there exists a supersolution $w_x$ in the domain $P(0,\nu_x)$ with $w_x(y+R\nu_x) \to \infty$ as $R \to \infty$.  In order for all of the continuity results to hold there will need to be a uniform lower bound on the growth of $w_x$ for $x \in \partial \Omega$.  For example this will hold for the above equation with the function $w_x(y) = y \cdot \nu_x$ as long as $\operatorname{div}_yA(x,y) \cdot \nu_x \geq 0$.
\end{rem}

\subsection{Outline of the Paper}
In Section \ref{sec: Set up} we discuss various notations and previous results which will be used throughout the paper.  For readers well versed in the field of homogenization of uniformly elliptic equations most of the material in Section \ref{sec: Set up} except for Section \ref{subsec: Rational and Irrational Directions} should be familiar.  In particular we will discuss in detail the properties of $\integer^n$ periodic functions restricted to hyperplanes with rational and irrational normal directions.  These results are essentially refinements of the classical equidistribution theorem of Weyl \cite{Weyl1910} on irrational rotations of the torus.  Then we describe precisely the assumptions on the differential operators we will consider and the associated comparison, regularity and interior homogenization results.

In Section \ref{sec: The Cell Problem} we analyze the cell problem \eqref{eqn: cell prob 1}.  We prove Theorem \ref{thm: main thm} part (i) and show properties of the homogenized boundary condition $\overline{\mu}$ in particular the continuity properties discussed above.  We show that the convergence of the solutions of the cell problem $v_{x_0,\nu}(y+R\nu)$ to $\overline{\mu}$ is sufficiently uniform at irrational directions to show the claimed continuity.  We also give an example which demonstrates the nonlinearity of the homogenized boundary condition when $F$ is nonlinear. 

In Section \ref{sec: Homogenization in General Domains} we use the results of the previous section to show the homogenization.  In particular we use the stability of the uniform convergence $v_{x_0,\nu}(y+R\nu) \to \overline{\mu}$ at irrational directions to show that the upper and lower half-relaxed limits of the $u^\e$ are respectively sub and supersolutions of \eqref{eqn: hom eqn 2}.  Then we use the result of Theorem \ref{thm: ABP type 1} to conclude that the upper and lower half-relaxed limits are the same and so the homogenization holds.

\subsection{Acknowledgements}  I would like to thank my advisor Inwon Kim for reading this article and providing many helpful comments.  I would also like to thank the anonymous referees for their suggestions and comments which helped to improve the clarity of the paper in various places.

\section{Set Up}
\label{sec: Set up}
\subsection{Notations} \label{subsec: Notations}
Let us make a note of the notational conventions we will follow in this work.  The numbers $0<\lambda<\Lambda$ will always refer to the ellipticity constants.  The meaning of the term ellipticity constants will be explained in Section \ref{subsec: Assumptions on the Operators}. The number $\alpha \in (0,1)$ will almost always refer to the H\"{o}lder continuity of whatever Dirichlet boundary data is under consideration.  Let $\Omega \subset \real^n$ the associated H\"{o}lder spaces are,
$$ C^\alpha(\Omega) := \left\{ \psi : \Omega \to \real : \|\psi\|_{C^\alpha(\Omega)}:=\sup_{x \in \Omega}|\psi(x)|+|\psi|_{C^\alpha(\Omega)}<+\infty\right\} $$
where $|\psi|_{C^\alpha}$ is the H\"{o}lder semi-norm defined by,
$$|\psi|_{C^\alpha}:= \sup_{x \neq y \in \Omega} \frac{|\psi(x)-\psi(y)|}{|x-y|^\alpha}.$$
We also define the oscillation of a real valued function $\psi$ on $\real^n$ over a set $E \subset \real^n$, 
$$\osc_E \psi :=\sup_E\psi - \inf_E \psi.$$
Constants denoted by $C$ and $c$ will denote universal constants, they depend only on $\lambda,\Lambda$ and the dimension $n$, they may change value from line to line and sometimes even within the same line.  Let $A$ and $B$ be two quantities, we write
$$ A \lesssim B \ \hbox{ to mean } \ A \leq C B \ \hbox{ where } C \hbox{ is a universal constant}. $$
If we want to emphasize the dependency of a constant on a quantity $A$ we will write $C(A)$.  
\subsection{Rational and Irrational Directions}\label{subsec: Rational and Irrational Directions}
 We define the following classes of normal directions with respect to the $\integer^n$ lattice,
\begin{DEF}
\begin{enumerate}[(i)]
\item $\nu \in S^{n-1}$ is called a rational direction if $\nu \in \real \integer^n$.
\item $\nu \in S^{n-1}$ is called an irrational direction if it is not a rational direction.
\end{enumerate}
\end{DEF}
Given a $\integer^n$ periodic function $\psi$ we we will be interested in the properties of $\left. \psi \right|_H$, the restriction of $\psi$ to the hyperplane $H = \{x \in \real^n: (x-x_0)\cdot \nu = 0\}$.  When $\nu$ is an irrational direction the distribution of the values taken by $\left. \psi \right|_H$ does not change too much under translation of $H$ (changing $x_0$).  This `fact' will be very important to understanding the homogenization of the boundary data in half spaces.

\medskip

 We will first give a heuristic description of the issue at hand.  Let us think of $\psi$ as a function on the unit periodicity cell $[0,1)^n$.  Then we may consider the parts of the unit periodicity cell which are cut by the hyperplane $H$,
$$H_{per} = \{ \tau \in [0,1)^n: \tau + z \in H \ \hbox{ for some } \ z \in \integer^n\} = \bigcup_{z \in \integer^n} (H+z) \cap [0,1)^n.$$ 
For hyperplanes $H$ whose normal direction is rational the union defining $H_{per}$ will be equivalent to a finite union.  As a result of this the restriction $\left. \psi \right|_H([0,1)^n)$ may be only a non-dense subset of all the values taken by $\psi$.  Moreover, this subset will be highly variable under changing $x_0$.  For example this issue is immediately evident in $\real^2$ when we consider the periodic function,
$$\psi(x) = \psi(x_1) = \sin 2\pi x_1 ,$$
and the hyperplanes $H = \{x_1 = a\}$.
\medskip

When the direction is irrational $H_{per}$ will be `uniformly dense' (in a sense to be made precise) in $[0,1)^n$.  This is exactly equivalent to the well known irrational rotation of the torus example in the case $n=2$.   In this case changing $x_0$ will not change the overall distribution of the values taken by $\left. \psi \right|_H$.  To be slightly more precise, this is the situation in which homogenization of the oscillating boundary condition $\psi(x/\e)$ will hold.

\medskip

Now we proceed to give a precise description of the above heuristics in a way that will be useful to us in the future.  Parts of the following discussion are borrowed from \cite{CKL12}.  The reason we repeat it here is we will need to make use of a more refined version of the results outlined there.
\begin{DEF} 
A bounded sequence of real numbers $(x_j)_{j=1}^\infty$ is said to be \textit{equidistributed} in an interval $[a,b]$ if for any $[c,d] \subseteq [a,b]$ we have,
$$\lim_{n \to \infty} \frac{|\{x_1,...,x_n\} \cap [c,d]|}{n} = \frac{d-c}{b-a}.$$
\end{DEF}
For $x \in \real$ let $[x]$ be the largest integer less than or equal to $x$.
\begin{DEF}
A sequence is \textit{equidistributed mod $1$} if the sequence $(x_j-[x_j])_{j \in \mathbb{N}}$ is equidistributed in $[0,1]$.
\end{DEF}
\begin{thm} \textup{(Weyl's equidistribution theorem \cite{Weyl1910})}
If $x$ is irrational then $(jx)_{j \in \mathbb{N}}$ is equidistributed mod $1$.
\end{thm}
In order to make quantitative estimates we introduce the notion of the discrepancy.  The following definition is from the book \cite{KN74} via \cite{CKL12}.  For a subset $E\subset [0,1]$, a natural number $N$ and a sequence $(x_j)_{j \in \mathbb{N}} \subset [0,1]$ call $A(E;N)$ the number of elements of the sequence $x_j$ for $1 \leq j \leq N$ contained in $E$. 
\begin{DEF}
 For any $x \in [0,1]$ let $x_j = jx-[jx]$ and define the discrepancy function,
$$D_N(x) = \sup_{E=[a,b) \subseteq [0,1]} \left| \frac{A(E;N) }{N}-|E|\right|.$$

\end{DEF}
Then by Weyl's equidistribution theorem $D_N(x) \to 0$ as $N \to \infty$ for each $x$ irrational.  Now let us replace $D_N(x)$ by a function equivalent up to constants which is continuous at every irrational.  We define the modified discrepancy function as in \cite{KN74},
$$D_N^* = D_N^*(x_1,...,x_N) = \sup_{0<a\leq 1}\left| \frac{A([0,a);N) }{N}-a\right|.$$
Then we have the following properties for $D_N^*(x)$,
\begin{lem}\label{lem: discrep prop} (Properties of $D_N^*$)
\begin{enumerate}[(i)]
\item (Theorem 1.3 in \cite{KN74}) The discrepancies $D_N$ and $D_N^*$ are equivalent up to constants, 
$$D_N^* \leq D_N \leq 2D_N^*.$$
\item (Theorem 1.4 in \cite{KN74}) Let $x_1 \leq x_2 \leq \cdot\cdot\cdot \leq x_N$ be in $[0,1)$.  Then their discrepancy $D_N^*$ is given by,
$$D_N^* = \frac{1}{2N}+\max_{i=1,...,N} \left|x_i - \frac{2i-1}{2N}\right|.$$
\item Define the modified discrepancy function for all $x \in \real$
$$D_N^*(x):= D_N^*(x-[x],2x-[2x],...,Nx-[Nx]), $$
then this function is continuous in a neighborhood of every irrational $x \in \real \setminus \rational$. 
\end{enumerate}
\end{lem}
\begin{proof}
The proofs of the first two parts can be found in the book \cite{KN74}.  We prove part (iii).  Given an irrational $x$ and an integer $N>0$ let us call
$$ \delta = \min \{|jx -kx| : 1\leq j<k\leq N\}>0.$$
Let $\sigma : \{1,...,N\} \to \{1,...,N\}$ be the permutation which orders the sequence $x-[x],...,Nx-[Nx]$, that is,
$$ \sigma(1)x -[\sigma(1)x]< \sigma(2)x-[\sigma(2)x] < \cdot \cdot \cdot < \sigma(N) x-[\sigma (N)x] .$$
Then for $|y-x|<\frac{\delta}{2N}$ we have that the above ordering is preserved for $y$, in fact for each $ 1 \leq j <N$,
\begin{align*} 
\sigma(j) y &\leq \sigma(j)x +\sigma(j)|y-x| < \sigma(j)x + \delta/2  \\
&\leq \sigma(j+1)x - \delta/2 < \sigma(j+1)x-\sigma(j+1)|y-x| \\
& \leq \sigma(j+1)y.
\end{align*}
Then from part (ii) we have that for $|y-x|<\frac{\delta}{2N}$ the discrepancy is has the form,
$$ D_N^*(y) = \frac{1}{2N}+\max_{i=1,...,N} \left|\sigma(i)y-[\sigma(i)y] - \frac{2i-1}{2N}\right|,$$
which, as a maximum of finitely many continuous functions which are continuous in a neighborhood of $x$, is also continuous in a neighborhood of $x$.
\end{proof}

\medskip

Now let $\nu = (\nu_1,...,\nu_n) \in S^{n-1}$ be a direction, and let $1 \leq i \leq n$ be the largest index corresponding to the largest component of $\nu$, that is,
$$ i = \max \{1\leq j \leq n : |\nu_j| = |\nu|_{\ell^\infty} \}.$$
Let $H_\nu = \{x \in \real^n : x \cdot \nu = 0\}$ be the hyperplane through the origin orthogonal to $\nu$.  For $1\leq j \leq n$ define $m_j(\nu)$ to be the slope of $H_\nu$ in the plane $\mathop{span}\{e_i,e_j\}$, 
\begin{equation}\label{eqn: ratios}
 0= -m_j(\nu)|\nu_i|+|\nu_j| = -m_j(\nu)|\nu|_{\ell^\infty}+|\nu_j|. 
 \end{equation}
The first equality indicates that one of the vectors $m_j(\nu)e_i\pm e_j$ is in the hyperplane $H_\nu$.  The second equality makes manifest that $m_j(\nu)$ thus defined are continuous in $\nu$ and $m_j(\nu) \in [0,1]$.  Now for every $\e>0$ and some $\gamma \in (0,1)$ we define,
\begin{equation}\label{eqn: omega def}
\omega_\nu(N):= 2\min_{1\leq j \leq n}D_{[N]}^*(m_j(\nu)) \ \hbox{ for } \ N>1.
\end{equation}
For irrational directions $\nu \in S^{n-1} \setminus \integer^n$ at least one of the $m_j(\nu)$ will be irrational.  As a result, $\omega_\nu(N) \to 0$ as $N \to \infty$ when $\nu$ is an irrational direction.  Furthermore, for $N$ sufficiently large, the $\min$ in \eqref{eqn: omega def} above will actually be the same as the $\min$ over $j$ such that $m_j(\nu) \in \real \setminus \rational$.  Then for each such $N>1$ fixed $\omega_\cdot (N)$ will be continuous in a $S^{n-1}$ neighborhood of $\nu$ with size depending on $N$ (it is a minimum of finitely many compositions of continuous functions).  

\medskip

Now we state a lemma which quantifies our heuristic argument that hyperplanes with irrational normals are uniformly dense modulo $\integer^n$ in $[0,1)^n$.  The lemma is from \cite{CKL12}, but we state it in a modified form which will be useful to us.  Due to the modifications we will also present the proof.

\begin{lem}\label{lem: irrational directions} \textup{(Lemma 2.7 in \cite{CKL12}) }  For $\nu \in S^{n-1}$ and $x_0 \in \real^n$ we define $H(x_0)  = \{ x \in \real^n: (x-x_0) \cdot \nu =0\}$.  Let $\omega_\nu: \integer_{>1} \to\real_+$ be as defined in \eqref{eqn: omega def}, then the following hold: 
\begin{enumerate}[(i)]
\item There exists a dimensional constant $C=C(n)>0$ such that the following is true: for any $x \in H(x_0)$ and any $N>1$ there exists $y \in \real^n$ such that
$$|x-y| \leq C(n)N \ \hbox{ with } \ y-x_0 \in \integer^n, $$
and
$$ \textup{dist}(y,H(x_0)) < \omega_\nu(N).$$
\item For any $\delta>0$, there exists $z \in \integer^n$ such that, 
$$\textup{dist}(z,H(x_0)) \leq \inf_{N>0} \omega_\nu(N) +\delta. $$
\end{enumerate}
\end{lem}
\begin{proof}
Let $\nu$ be any direction in $S^{n-1}$ and $N>1$.  Without loss we may translate so that $x = 0$ and assume that $\nu_n = |\nu|_{\ell^\infty} \geq 1/\sqrt{n}$.  Let $j$ such that $\omega_\nu(N) = D_{[N]}^*(m_j(\nu))$.  Recall that one of the vectors $m_j(\nu)e_n\pm e_j$ is in the hyperplane $H = H(0)$, let us assume it is $m_j(\nu)e_n+ e_j$, the other case will work out the same.  Now from the Lemma \ref{lem: discrep prop}  we know that the the discrepancy $D_{[N]}(m_j(\nu)) \leq \omega_\nu(N)$, and therefore every interval $[a,b) \subset [0,1)$ with length at least $\omega_\nu(N)$ contains at least one of the $km_j(\nu)-[km_j(\nu)]$ with $k \leq N$.

\medskip

Let $s = s_ne_n$, with $s_n \in [0,1)$, then by the above argument there exists $1\leq k \leq N$ such that $|km_j(\nu)-[km_j(\nu)] -s_n| \leq \omega_\nu(N)$, and therefore,
$$ |k(m_j(\nu)e_n+e_j)-s| \bmod \integer^n \leq \omega_\nu(N). $$
In particular there exists $y \in H_\nu$ such that $|y-s| \bmod \integer^n \leq \omega_\nu(N)$ and $|y| \leq \sqrt{2}N$. 

\medskip

Now let us consider $s \in [0,1)^n$ arbitrary.  Then let us take 
$$ s' = \frac{s \cdot \nu}{\nu_n}e_n  \ \hbox{ so that } \ s-s' \in H. $$ 
  Now, by the above arguments, there exists $y' \in H$ such that $|y'-s'| \bmod \integer^n \leq \omega_\nu(N)$ and $|y'| \leq \sqrt{2}N$.  Then let $y = y'+(s-s')$ which is still in $H$ and $|y| \leq |y'|+|s-s'| \leq 2\sqrt{n} N$ and, moreover,
$$ |y-s| \bmod \integer^n = |y'-s'| \bmod \integer^n \leq  \omega_\nu(N).$$

Now let $s = x_0 \bmod \integer^n$, then there exists $\tilde{y} \in H$ as above with 
$$|\tilde{y}| \leq 2\sqrt{n} N \ \hbox{ and } \ |\tilde{y} - x_0| \bmod \integer^n =  |\tilde{y} - s| \bmod \integer^n \leq \omega_\nu(N).$$
Now let $y = \tilde{y}+(s-\tilde{y} \bmod \integer^n)$, then we have $|y| \leq (2\sqrt{n}+1)N$ and
$$ y \bmod \integer^n = s = x_0 \bmod \integer^n \ \hbox{ and } \ \text{dist}(y,H) \leq \omega_\nu(N). $$
This completes the proof of part (i).

\medskip

The proof of part (ii) is very similar.  Let $\delta>0$ and $N$ sufficiently large so that $\omega_\nu(N) \leq \inf_{N>0} \omega_\nu(N) +\delta$.  Let $s = -x_0$ and by the arguments above there exists $1 \leq k \leq N$ such that
$$ |k(m_j(\nu)e_n+e_j)+x_0| \bmod \integer^n \leq \omega_\nu(N) \leq \inf_{N>0} \omega_\nu(N) +\delta. $$
Since $x_0+k(m_j(\nu)e_n+e_j)$ is in $H(x_0)$ this completes the proof of (ii).
\end{proof}

\subsection{Assumptions on the Operators}\label{subsec: Assumptions on the Operators}
Now we give the technical assumptions on the differential operators under which we can solve the cell problem.  The operator $F$ which arises in the cell problem is obtained as a scaling limit of the general operators $F_\e$ being homogenized over.  This connection will be made more explicit later in Lemma \ref{lem: scale one op}.    We will work in the class of fully nonlinear uniformly elliptic operators without gradient dependence.  This is mostly for convenience and we believe that our framework can be extended to the gradient dependent case by using the arguments of \cite{BarlesMironescu12}.  

\medskip

Let $\mathcal{M}^n$ be the class of $n \times n$ symmetric matrices with real entries.  Then the differential operator $F : \mathcal{M}^n \times \real^n  \to \real$ will be assumed to satisfy the following properties:
\begin{enumerate}[(F1)]
\item (Lipschitz continuity)  $F$ is locally Lipschitz continuous in $\mathcal{M}^n \times \real^n$, and moreover there exists $C>0$ so that for all $x,y \in \real^n$ and $M,N \in \mathcal{M}^n$,
$$ |F(M,z)-F(N,y)| \leq C(|z-y| (1+\|M\|+\|N\|)+\|M-N\|).$$
\item (Uniform Ellipticity) There exists $\Lambda>\lambda>0$ such that for all $y \in \real^n$,
$$ \lambda \Tr( N) \leq F(M,y)-F(M+N,y) \leq \Lambda \text{Tr}( N) \ \hbox{ for all } M,N \in \mathcal{M}^n \ \hbox{ such that } N \geq 0.$$
\item (Periodicity) For every $z \in \integer^n$, $y \in \real^n$ and $M \in \mathcal{M}^n$ we have,
$$F(M,y+z) = F(M,y).$$
\item (Homogeneity) For all $t >0$, $M \in \mathcal{M}^n$ and $y \in \real^n$ we have $F(tM,y) = tF(M,y)$ and, in particular $F(0,y)=0$. 
\end{enumerate}

A wide class of operators of this form come from the theory of optimal control of diffusion processes, for example the Hamilton-Jacobi-Bellman operators,
$$ F(M,y) = -\sup_{\alpha \in \mathcal{A}} \Tr(A^\alpha(y)M) \ \hbox{ where } \ \lambda  \leq A^\alpha(y) \leq \Lambda \ \forall \ \alpha \in \mathcal{A}. $$
In much greater generality there are Isaacs operators coming from differential games which have the form,
$$ F(M,y) = -\inf_{\beta \in \mathcal{B}}\sup_{\alpha \in \mathcal{A}} \Tr(A^{\alpha\beta}(y)M) \ \hbox{ where } \ \lambda  \leq A^{\alpha\beta}(y) \leq \Lambda \ \forall \  (\alpha,\beta) \in \mathcal{A}\times \mathcal{B}.$$
In fact, all operators satisfying (F1)-(F4) can be written as an Isaacs operator, the proof is a basic exercise.

\medskip

The operator $F_\e$ is assumed to satisfy (F1) with the local Lipschitz continuity in both the $x$ and $y$ variables, the uniform ellipticity (F2) and (F4) the $F_\e$ are $\integer^n$ periodic in $y$.  The most general form of such operators is,
\begin{equation}
 F_\e(M,x,y) = f(x,y)-\inf_{\beta \in \mathcal{B}}\sup_{\alpha \in \mathcal{A}} \Tr(A^{\alpha\beta}(x,y)M) \ \hbox{ where } \ \lambda  \leq A^{\alpha\beta} \leq \Lambda \ \forall \  (\alpha,\beta) \in \mathcal{A}\times \mathcal{B}.
 \end{equation}
Here the $A^{\alpha,\beta}$ have uniformly bounded Lipschitz norm, are periodic in their $y$ variable, and satisfy the uniform ellipticity condition
$$\lambda  \leq A^{\alpha\beta} \leq \Lambda \ \forall \  (\alpha,\beta) \in \mathcal{A}\times \mathcal{B}.$$
The function $f$ is in  $C^{0,1}(\real^n \times \real^n)$ and is $\integer^n$ periodic in its second variable.

The standard notion of weak solution for the equation $F(D^2u,y)=0$ in a domain $\Omega \subseteq \real^n$ is the viscosity solution.  The usual reference on the theory of viscosity solutions is \cite{CIL92}.  
\begin{DEF} (Viscosity Solutions) 
\begin{enumerate}[(i)]
\item $u: \Omega \to \real$ upper semi-continuous is a subsolution of $F(D^2u,y) \leq 0 $ if for any $\phi \in C^2(\Omega)$ such that $u-\phi$ has a local max at $y_0 \in \Omega$,
$$ F(D^2\phi(y_0),y_0) \leq 0.$$
\item $u: \Omega \to \real$ lower semi-continuous is a supersolution of $F(D^2u,y) \geq 0 $ if for any $\phi \in C^2(\Omega)$ such that $u-\phi$ has a local min at $y_0 \in \Omega$,
$$ F(D^2\phi(y_0),y_0) \geq 0.$$
\item $u: \Omega \to \real$ continuous is a solution of $F(D^2u,y) = 0$ if it is both a subsolution and a supersolution.
\end{enumerate}
\end{DEF}

The Dirichlet boundary data can be defined in the classical sense for viscosity solutions.  We say that $u$ is a subsolution of $F(D^2u,y) = 0$ in $\Omega$ with Dirichlet data $\psi : \partial\Omega \to \real$ if $u$ is a subsolution and
$$ \limsup_{z \to y} u(z) \leq \psi(y) \ \hbox{ for all } \ y \in \partial\Omega. $$
Supersolutions are defined analogously and solutions satisfy both the sub solution and the super solution property.  

\medskip

One very useful property of viscosity solutions is their stability under uniform convergence.  Let $u_j : \overline{\Omega} \to \real$ with $j \in \mathbb{N}$ for some domain (possible unbounded) $\Omega \subset \real^n$.  Let us defined the upper and lower half-relaxed limit operations,
$$ \left.\limsup\right.^* u_j(x) = \lim_{\e \to 0} \ \sup\{ u_j(y) : y \in \overline{\Omega} \ \hbox{ and } \ |y-x|+j^{-1} \leq \e \}$$
and
$$ \left.\liminf\right.^* u_j(x) = \lim_{\e \to 0} \ \inf\{ u_j(y) : y \in \overline{\Omega} \ \hbox{ and } \ |y-x|+j^{-1} \leq \e \}.$$
One can also have $u_j : \Omega_j \to \real$ with the $\Omega_j \to \Omega$ in the Hausdorff topology on compact subsets of $\real^n$.  One can easily check that if the upper and lower half-relaxed limits agree then the $u^j$ converge locally uniformly to their common value.  In any case we have the following stability result for viscosity sub and supersolutions under the half-relaxed limit operations, see \cite{CIL92} for more details,
\begin{lem}\label{lem: uniform stability} \textup{(Stability of viscosity solutions)}
Let $F_j(M,y): \mathcal{M}^n \times \real^n \to \real$ be a sequence of operators satisfying (F1) and (F2) which converge to an operator $F(M,y)$ locally uniformly in $M$ and uniformly in $y \in \real^n$.  Let $\Omega$ be a domain in $\real^n$ and let $u_j: \Omega \to \real$ be a bounded above (below) sequence of subsolutions (supersolutions) of,
$$ F_j(D^2u_j,y) \leq 0 \ \hbox{ in } \ \Omega  \ \ \ ( \ F_j(D^2u_j,y) \geq 0 \ \hbox{ in } \ \Omega \ ). $$
Define $u^*$ the upper half-relaxed limit of the $u^j$ (resp. $u_*$ the lower half-relaxed limit), then we have,
$$ F(D^2u^*,y) \leq 0 \ \hbox{ in } \ \Omega  \ \ \ ( \ F(D^2u_*,y) \geq 0 \ \hbox{ in } \ \Omega \ ). $$
\end{lem}
Then we state the interior homogenization result in terms of the half-relaxed limits:
\begin{thm} (Interior homogenization \cite{Evans92,Evans89,CSW05})\label{thm: interior hom}
Let $F_\e(M,x,y)$ satisfying (F1) with the local Lipschitz continuity in both the $x$ and $y$ variables, (F2) and $F_\e$ are $\integer^n$ periodic in $y$.  Let $\Omega$ a bounded domain in $\real^n$ and $u^\e$ satisfying,
\begin{equation}
\left\{
\begin{array}{lll}
F_\e(D^2u^\e,x,\tfrac{x}{\e}) \leq 0  & \hbox{in} & \Omega \\ \\
u^\e(x) \leq M & \hbox{on} & \partial\Omega
\end{array}\right.
\end{equation}
for some $M>0$.  Then there exists an operator $\overline{F}(M,x)$ such that $u^* = \left.\limsup\right.^* u^\e $ is a subsolution of,
\begin{equation}
\left\{
\begin{array}{lll}
\overline{F}(D^2u^*,x) \leq 0  & \hbox{in} & \Omega \\ \\
u^*(x) \leq M & \hbox{on} & \partial\Omega.
\end{array}\right.
\end{equation}
The analogous result holds for supersolutions as well.
\end{thm}
\medskip

\subsection{Maximal Operators and Comparison}\label{subsec: Maximal Operators and Comparison}
Now we recall the Pucci maximal operators associated with the class of uniformly elliptic operators with constants $\lambda, \Lambda$.  The basic results about the Pucci maximal operators can be found in the book \cite{CaffarelliCabre95}.  For $M \in \mathcal{M}^n$ we can always decompose $M=M_+-M_-$ with $M_\pm \geq 0$ and $M_+M_-=0$.  We then define,
\begin{equation}\label{eqn: maximal operators}
 \mathcal{P}^+(M) = \Lambda \Tr M_+-\lambda \Tr M_- \ \hbox{ and } \ \mathcal{P}^-(M) = \lambda \Tr M_+-\Lambda \Tr M_- .
 \end{equation}
 Then from (F2) we may derive for $M,N\in \mathcal{M}^n$ and $x \in \real^n$,
 \begin{equation}\label{eqn: ellipticity w/ pucci}
 -\mathcal{P}^+(M-N)\leq F(M,y)-F(N,y) \leq -\mathcal{P}^-(M-N).
 \end{equation}
 Note that by the inequality \eqref{eqn: ellipticity w/ pucci} any viscosity solution of $F(D^2u,x)=0$ will satisfy, in the viscosity sense,
$$ -\mathcal{P}^+(D^2u) \leq 0 \ \hbox{ and } \ -\mathcal{P}^-(D^2u) \geq 0.$$
Moreover, \eqref{eqn: ellipticity w/ pucci} implies that the difference $w=u-v$ of a classical supersolution $u$ and a classical subsolution $v$ is itself a supersolution of the equation,
 $$ -\mathcal{P}^-(D^2w) \geq 0. $$
 This is also true of only semicontinuous viscosity sub and supersolutions using the method of sup and inf-convolutions originally used for this purpose by Jensen \cite{Jensen88}.  We state the Theorem in a slightly different form closer to what appears in \cite{CaffarelliCabre95},
 \begin{thm}\textup{(Comparison)} \label{thm: comparison pucci} Let $\Omega$ be a domain in $\real^n$ not necessarily bounded.  Let $F$ satisfy (F1) and (F2).  Suppose that $u$ and $v$ satisfy, in the viscosity sense,
$$F(D^2u,y) \leq F(D^2v,y)  \ \hbox{ in } \ \Omega. $$
Then $w=u-v$ is a subsolution of,
$$ -\mathcal{P}^+(D^2w) \leq 0 \ \hbox{ in } \ \Omega. $$
\end{thm}
In bounded domains the comparison principle for the Dirichlet problem follows from Theorem \ref{thm: comparison pucci} in a standard way.  A comparison principle for solutions of uniformly elliptic equations in half spaces with sublinear growth at infinity then follows from the following localization lemma which will be useful to us later,
\begin{lem}\label{lem: localization} \textup{(Localization)} Let $\delta,\e>0$ and $L>1$ and suppose that $v$ is a viscosity solution of
\begin{equation}\label{eqn: localization}
\left\{
\begin{array}{lll}
-\mathcal{P}^+(D^2v) \leq 0  & \hbox{ in } & P(0,e_n) \cap Q_L\\ \\
v(y) \leq \delta & \hbox{ on } & \partial P(0,e_n) \cap Q_L \\ \\
v(y) \leq L^{1-\e} & \hbox{ on } & \partial Q_L \cap \overline{P},
\end{array}\right.
\end{equation}
where $Q_L$ is the cylindrical region,
$$Q_L = \{|x'| \leq L\} \times \{ 0 \leq x_n \leq L\} \ \hbox{ with } \ x' = x-x_ne_n.$$
  Then we have,
$$ v(y) \leq  \delta+2\tfrac{\Lambda}{\lambda}nL^{-\e} \ \hbox{ in } \ Q_1.$$
\end{lem}
\begin{rem}
Since the Pucci maximal operators are translation and rotation invariant the same result holds for the analogous problem in the half space $P(x_0,\nu)$ for any $\nu \in S^{n-1}$ and $x_0 \in \real^n$.
\end{rem}
\begin{proof}
  Consider the following barrier $\phi$,
$$ \phi(x) = L^{-1-\e}(|x'|^2-2\tfrac{\Lambda}{\lambda}(n-1)x_n^2)+\delta+(2\tfrac{\Lambda}{\lambda}(n-1)+1)L^{-\e}x_n.$$
It is straightforward to check that $\phi(x) \geq v(x)$ on $\partial Q_R$ and that $\phi$ is a smooth supersolution of $-\mathcal{P}^+(D^2\phi) \geq 2 \Lambda(n-1)L^{-1-\e} $ in $x_n>0$.  From the definition of viscosity solution $v \leq \phi$ in $Q_L$.  Therefore since
$$ \phi(x) \leq \delta+2(\tfrac{\Lambda}{\lambda}(n-1)+1)L^{-\e} \ \hbox{ in } \ Q_1,$$
we get the result.

\end{proof}

A straightforward application of the localization lemma is the following comparison principle for bounded solutions of the half-space problem \eqref{eqn: cell prob 1}:

\begin{lem}\label{lem: comparison}
Let $\psi_1,\psi_2 \in C^\alpha(\real^n)$. Let $v_1$ be a bounded upper semi-continuous subsolution of \eqref{eqn: cell prob 1} with Dirichlet data $\psi_1$, and $v_2$ be a bounded lower semi-continuous supersolution of \eqref{eqn: cell prob 1} with Dirichlet data $\psi_2$.  Then,
$$ \sup_{y \in P} \ (v_1(y)-v_2(y))_+ \leq \sup_{y \in \partial P} (\psi_1(y)-\psi_2(y))_+.$$
\end{lem}

\subsection{Regularity Results}
The Pucci maximal operators govern the worst possible behavior for solutions of equations satisfying the ellipticity condition (F2).  In addition to the role they play in the comparison result, Theorem \ref{thm: comparison pucci}, they are useful because regularity results hold uniformly in the ellipticity class $(\lambda,\Lambda)$.  In particular we have the following result, for more details see the book \cite{CaffarelliCabre95}:
\begin{lem}\label{lem: int estimates}
\textup{(Interior H\"{o}lder estimates)}  Assume (F2) and (F4).  Let $x_0 \in \real^n$ and $v$ be a continuous viscosity solution of 
$$ -\mathcal{P}^+(D^2v) \leq 0 \ \hbox{ and } \ -\mathcal{P}^-(D^2v) \geq 0.$$
 in $B_r(x_0)$. For every $\alpha \in (0,1)$ there exists $C = C(\lambda,\Lambda,n,\alpha)>0$ such that,
$$ \sup_{x,y \in B_{r/2}(x_0)} \frac{|v(x)-v(y)|}{|x-y|^\alpha} \leq Cr^{-\alpha} \sup_{x \in B_r(x_0)} v(x).$$
\end{lem}

By using barriers one can prove the following H\"{o}lder estimate up to the boundary,
\begin{lem}\label{lem: bdry estimates}
\textup{(Boundary H\"{o}lder estimates)}  
Let $g \in C^{\alpha}(\real^n)$ for some $\alpha \in (0,1)$ and suppose that $u$ is a viscosity solution of
\begin{equation}
\left\{
\begin{array}{lll}
F(D^2u) = 0 &\hbox{ in } & B_1^+ = \{|x| < 1, \ x_n>0\} \\ \\
u = g(x) & \hbox{ on } & T=\partial B_1^+ \cap \{x_n = 0\}.
\end{array}
\right.
\end{equation}
Then we have the following estimate up to the boundary,
$$\|u\|_{C^{\alpha}(\overline{B^+_{1/2}})} \leq C(n,\lambda,\Lambda)(\|g\|_{C^{\alpha}(T)}+\sup_{B_1^+} |u|).$$
\end{lem}

\section{The Cell Problem}\label{sec: The Cell Problem}

Now we begin discussing the actual homogenization problem.  First we will analyze the cell problem and identify the homogenized boundary condition.  We will work in the following half-space type domains in $\real^N$,
$$P(\nu,x_0) = \{ x \in \real^n : (x-x_0) \cdot \nu >0 \}.$$  
Let $\psi : \real^n \to \real$ be $C^\alpha(\real^n)$ for some $\alpha \in (0,1)$ and periodic with respect to the $\integer^n$ lattice.  For convenience in the future we will sometimes refer to this class of functions as $C^{\alpha}_{per}(\real^n)$.  From here on $\alpha$ will always refer to the exponent of H\"{o}lder continuity of the Dirichlet boundary data for the problem under consideration.  Now we consider solving the following Dirichlet problem in $P(\nu,x_0)$, which we will call the cell problem,
\begin{equation}\label{eqn: cell prob}
\left\{
\begin{array}{lll}
F(D^2v_{x_0,\nu},y) = 0  & \hbox{ in } & P(\nu,x_0) \\  \\
v_{x_0,\nu}(y) = \psi(y) & \hbox{ on } & \partial P(\nu,x_0).
\end{array}\right.
\end{equation}
A standard application of Perron's method in combination with the Localization Lemma \ref{lem: localization} and the Comparison Principle \ref{lem: comparison} give a unique bounded solution to the cell problem.  Standard barrier arguments imply that the solution given by Perron's method achieves the boundary data continuously.  

\medskip

Now let us define the upper and lower homogenized boundary conditions arising from the cell problems \eqref{eqn: cell prob}, 
\begin{equation}
\begin{array}{l}
\mu^*(\nu,F,\psi) := \sup \{ \limsup_{R \to \infty} v_{x_0,\nu}(y+R\nu) : y \in \partial P(\nu,x_0), \ x_0 \in \real^n\} \\
\mu_*(\nu,F,\psi) := \inf \{ \liminf_{R \to \infty} v_{x_0,\nu}(y+R\nu) : y \in \partial P(\nu,x_0), \ x_0 \in \real^n\}. 
\end{array}
\end{equation}
Then we call $\overline{\mu} = \frac{1}{2}(\mu^*+\mu_*)$ to be the average of $\mu^*$ and $\mu_*$.  Notice that $\mu^*$, $\mu_*$ and $\overline{\mu}$ are invariant under translations of $\psi$.  In fact we could have equivalently defined them as envelopes over $\psi(\cdot+\tau)$ for $\tau \in [0,1)^n$ rather than over the location of the hyperplane $\partial P(\nu,x_0)$.  From now on we will work mostly with the quantity $\overline{\mu}(\psi,F,\nu)$ but the upper and lower homogenized boundary conditions will be mentioned occasionally as well since we will see that they are relevant at the rational directions.  

Now we will show that  $\overline{\mu}=\mu^*=\mu_*$ when $\nu$ is an irrational direction and that the convergence of $v_{x_0,\nu}(y+R\nu)$ to $\overline{\mu}$ is uniform in $x_0$ and $y$ as $R \to \infty$.  Furthermore we provide a precise estimate on the difference $\mu^*-\mu_*$ when $\nu$ is a rational direction in terms of the discrepancy in the form of $\omega_\nu$.  This lemma is analogous to Theorem 2.2 in \cite{BarlesMironescu12}.    
\begin{lem}  \label{lem: ID BC}\textup{(Identification of the homogenized boundary condition)} Suppose that $v_{x_0}=v_{x_0,\nu}$ is the unique bounded solution of \eqref{eqn: cell prob} set in $P(\nu,x_0)$ where $\nu$ is any direction in $ S^{n-1}$ and $x_0 \in \real^n$.  Let $\overline{\mu}(\psi,F,\nu)$ defined as above then we have for all $N>1$,
\begin{equation}\label{eqn: cell prob unif conv}
 \sup_{x_0 \in \real^n}\sup_{y \in \partial P} |v_{x_0}(y+R\nu) -\overline{\mu}| \leq C(n,\lambda,\Lambda)(|\psi|_{C^{\alpha}}+\osc \psi)((N/R)^{\alpha}+\omega_\nu(N)^\alpha) \ \hbox{ for } \ R \geq C N.\end{equation}
In particular, when $\nu$ is an irrational direction $\nu \in S^{n-1} \setminus \real\integer^n$ we have that $\mu^*=\mu_* = \overline{\mu}$ and for all directions $\nu \in S^{n-1}$,
\begin{equation}\label{eqn: mus diff}
\mu^*(\nu,F,\psi)-\mu_*(\nu,F,\psi) \leq C(n,\lambda,\Lambda)(|\psi|_{C^{\alpha}}+\osc \psi) \inf_{N>1} \omega_\nu(N)^\alpha.
\end{equation}
\end{lem}
\begin{rem}
The fact that the limit of $v_{x_0,\nu}(y+R\nu)$ as $R \to \infty$ does not depend on the location of the hyperplane $\partial P$ (i.e. on $x_0$) for irrational directions is essential for homogenization to hold.  As mentioned before this is heuristically because restrictions of $\psi$ to irrational hyperplanes see all the values of $\psi$.  The H\"{o}lder regularity of $\psi$ is not essential, but in general continuity seems to be necessary.
\end{rem}
\begin{proof}
We may assume without loss of generality due to invariance of the equation under adding a constant and scalar multiplication that $|\psi|_{C^\alpha} \leq 1$ and that $\|\psi\|_\infty = \tfrac{1}{2}\osc \psi \leq 1$.

1. First we show that $v_{0,\nu}(y+R\nu)$ converges along subsequences locally uniformly in $\real^n$ to constants.  For now let us call $v = v_{0,\nu}$, $P = P(\nu,0)$.  By the comparison principle we have the uniform bounds,
$$ -1 \leq v \leq 1. $$
Let $R_k \to \infty$ be a subsequence such that $v(R_k\nu)$ converges as $R_k \to \infty$ to some constant $\tilde{\mu}$.

\medskip

2.  Due to Lemma \ref{lem: irrational directions} for any $p = p'+R\nu \in \partial P+R\nu$ there exists $z \in \integer^n$ such that,
$$ |z-p'| \leq CN \ \hbox{ and } \ h:=\text{dist}(z,\partial P) \leq \omega_\nu(N).$$
  Note that we have $z\cdot \nu=\pm h$, we assume the $+$ sign, the proof will be very similar in the other case.  We define the translation of $v$,
$$\widetilde{v}(y) = v(y-z),$$
which is a viscosity solution of,
\begin{equation*}
\left\{
\begin{array}{lll}
 F(D^2\widetilde{v},y)= F(D^2\widetilde{v},y-z) = 0  & \hbox{ in } & P+h\nu \\ 
\widetilde{v}(y) = \psi(y-z) = \psi(y) & \hbox{ on } & \partial P+h\nu.
\end{array}\right.
\end{equation*}
  Now for $y \in \partial P+h\nu \subset \overline{P}$ we have,
$$|v(y)-\psi(y)| \leq |v(y)-v(y-h\nu)|+|\psi(y-h\nu)-\psi(y)| \leq C h^\alpha.$$
Therefore by comparison principle,
$$ \sup _{y \in P+h\nu}|v(y)-v(y-z)| \leq Ch^\alpha, $$
and then in particular, by the interior $C^\alpha$ estimates,
\begin{align*}
 |v(p+R\nu)-v(R\nu)|&\leq |v(p+R\nu)-v(z+R\nu)|+|v(z+R\nu)-v(R\nu)| \\
&\leq   \osc_{|y-p|\leq CN} v(y)+C\omega_\nu(N)^{\alpha} \leq C( \omega_\nu(N)^\alpha+(N/R)^\alpha),
\end{align*}
where $C$ is of course independent of the point $p$.  Then using comparison, 
$$ \osc_{y \in P} v(y+R\nu) \leq \osc_{y \in \partial P} v(y+R\nu) \leq C( \omega_\nu(N)^\alpha+(N/R)^\alpha).$$
Now for $R_k>R$, since $R_k\nu \in P+R\nu$,
$$\sup_{y \in P} |v(y+R\nu)-v(y+R_k\nu)| \leq C( \omega_\nu(N)^\alpha+(N/R)^\alpha).$$
Letting $R_k\to \infty$ in the above inequality we derive,
\begin{equation}\label{eqn: in lem est 1}
\sup_{y \in P(0,\nu)}|v_{0,\nu}(y+R\nu)-\tilde{\mu}| \leq C(\omega_\nu(N)^{\alpha} +(N/R)^{\alpha}).
\end{equation}

\medskip

3. Now we show that the estimate \eqref{eqn: in lem est 1} is independent of the spatial location of $\partial P$.  Let $x_0 \in \real^n$ and $v_{x_0}$ be the corresponding solution of \eqref{eqn: cell prob} in $P(\nu,x_0)$, and $v$ as before be the solution of \eqref{eqn: cell prob} in $P(\nu,0)$.  Let $N>1$, due to Lemma \ref{lem: irrational directions} there exists $z \in \integer^n$ such that
$$ \text{dist}(\partial P(\nu,x_0)+z,\partial P(\nu,0)) = h \leq \omega_\nu(N). $$
 Then $\widetilde{v}(y)=v(y-z)$ solves the following problem in $P(\nu,z)$,
\begin{equation*}
\left\{
\begin{array}{lll}
 F(D^2\widetilde{v},y)= F(D^2\widetilde{v},y-z) = 0  & \hbox{ in } & P(\nu,z) \\ 
\widetilde{v}(y) = \psi(y-z) = \psi(y) & \hbox{ on } & \partial P(\nu,z).
\end{array}\right.
\end{equation*}
Note that $\widetilde{v}$ trivially satisfies the estimate \eqref{eqn: in lem est 1}.  Let us assume that $z \in P(\nu,x_0)$, the other case will follow by a similar argument.  For each $y \in \partial P(\nu,z)$, we have that $y - h \nu \in \partial P(\nu,x_0)$, so by the $C^{\alpha}$ regularity up to the boundary for $v_{x_0}$ we get,
$$|v_{x_0}(y)-\psi(y)| \leq |v_{x_0}(y)-v_{x_0}(y - h\nu)|+|\psi(y - h\nu)-\psi(y)| \leq C h^ \alpha.$$
Therefore by the comparison principle Lemma \ref{lem: comparison},
$$ \sup_{y \in P(\nu,z)} |\widetilde{v}(y)-v_{x_0}(y)| \leq C \omega_\nu (N)^\alpha, $$
Combining the above estimate with estimate \eqref{eqn: in lem est 1} for $\widetilde{v}$, we get
\begin{equation}\label{eqn: in lem est 2}
\sup_{y \in P(0,\nu)}|v_{x_0,\nu}(y+R\nu)-\tilde{\mu}| \leq C(\omega_\nu(N)^{\alpha} +(N/R)^{\alpha}),
\end{equation}
where the constant $C$ depends only on the uniform ellipticity of $F$ by way of the interior $C^{\alpha}$ estimates, and in particular not on the arbitrary point $x_0$ or $N$.  Taking the supremum over $x_0 \in \real^n$ and we get the desired estimate except with $\tilde{\mu}$.  From here we simply note that $\tilde{\mu} \in [\mu_*,\mu^*]$ so estimate \eqref{eqn: cell prob unif conv} with $\tilde{\mu}$ replacing $\overline{\mu}$ still implies \eqref{eqn: mus diff} which in turn implies \eqref{eqn: cell prob unif conv} with $\overline{\mu}$ (possibly at the cost of increasing the universal constant $C(n,\lambda,\Lambda)$).

\end{proof}
We now list some properties of the homogenized boundary condition which can be derived in a straightforward way from the identification Lemma \ref{lem: ID BC}. 
\begin{lem}\label{lem: mu properties}\textup{(Properties of $\overline{\mu}$)}  The following properties all hold for $\mu^*$ and $\mu_*$ as well unless otherwise specified.  Recall that the distinction between the envelopes is only relevant at rational directions.
\begin{enumerate}[(i)]
\item (Continuity with respect to $F$) Fix $\psi \in C^{\alpha}_{per}(\real^n)$, $\nu \in S^{n-1}$ and let $F_j(M,y)$ be a sequence of operators satisfying (F1)-(F4) which converge to an operator $F_\infty(M,y)$ locally uniformly in $M$ and uniformly in $y$, then,
$$ \limsup_{j \to \infty} |\overline{\mu}(\psi,F_j,\nu) -\overline{\mu}(\psi,F_\infty,\nu)| \leq C(n,\lambda,\Lambda)(|\psi|_{C^{\alpha}}+\osc \psi)\inf_{N>0}\omega_\nu(N)^\alpha.$$
In particular this implies the continuity of $\overline{\mu}$ with respect to locally uniform convergence of its $F$ argument at all irrational directions $\nu \in S^{n-1} \setminus \real \integer^n$.
\item (Monotonicity) For any $\psi_1, \psi_2 \in C^{\alpha}_{per}(\real^n)$ and $F_1, F_2$ satisfying (F1)-(F4) we have,
$$ \psi_1 \leq \psi_2 \ \hbox{ and } \ F_2 \leq F_1 \  \hbox{ implies }  \ \overline{\mu}(\psi_1,F_1,\nu) \leq \overline{\mu}(\psi_2,F_2,\nu).$$
\end{enumerate}
For the following we fix $\nu \in S^{n-1}$ and an operator $F$ satisfying (F1)-(F4), and we call $\overline{\mu}(\cdot,F,\nu) = \overline{\mu}(\cdot)$ for simplicity of notation.
\begin{enumerate}[(i)]
\setcounter{enumi}{2}
\item (Continuity with respect to $\psi$)  For $\psi_1, \psi_2 \in C^\alpha_{per}(\real^n)$,
$$ |\overline{\mu}(\psi_1)-\overline{\mu}(\psi_1)| \leq \| \psi_1-\psi_2\|_{L^\infty}.$$
\item (Homogeneity) For any $t>0$ and $\psi \in C^{\alpha}_{per}(\real^n)$,
$$ \overline{\mu}(t\psi) = t \overline{\mu}(\psi) .$$
\item (Constants)  For $\psi \in C^{\alpha}_{per}(\real^n)$ and $c \in \real$,
$$\overline{\mu}(\psi+c) = \overline{\mu}(\psi)+c .$$
\item (Translation invariance) For any $\psi \in C^{\alpha}_{per}(\real^n)$ and any $\tau \in [0,1)^n$,
$$ \overline{\mu}(\psi(\cdot+\tau)) = \overline{\mu}(\psi).$$
\item (Sub/Super-additivity) If the operator $F$ is convex (concave) then, for any $\psi_1, \psi_2 \in C^{\alpha}_{per}(\real^n)$,
$$ \mu_*(\psi_1+\psi_2) \geq  \mu_*(\psi_1)+\mu_*(\psi_2) \ \ ( \hbox{ or } \ {\mu}^*(\psi_1+\psi_2) \leq  {\mu}^*(\psi_1)+{\mu}^*(\psi_2) ).$$
In particular, if $F$ is linear then $\overline{\mu}$ is linear as well at irrational directions.
\end{enumerate}
\end{lem}
\begin{proof}
All except (i) are easy applications of Lemma \ref{lem: ID BC} and the comparison principle.  Let us therefore consider the situation in (i).  As before we assume that $|\psi|_{C^\alpha} \leq 1$ and that $\|\psi\|_\infty = \tfrac{1}{2}\osc \psi $.  Let $\overline{\mu}_j = \overline{\mu}_j(\psi,F_j,\nu)$ and let $v_j$ for $j \in \mathbb{N}\cup\{\infty\}$ respectively solve the cell problems \eqref{eqn: cell prob} for $F_j$ in the domain $P=P(0,\nu)$.  By the method of half-relaxed limits \cite{BarlesPerthame88} the $v_j$ converge uniformly on compact subsets of $P \cup \partial P$ to $v_\infty$.  Let $\delta>0$ and let $N$ sufficiently large so that $\omega_\nu(N)^\alpha \leq c\delta+ \inf_{N>0} \omega_\nu(N)^\alpha$ and then taking $R\geq \delta^{-1/\alpha}CN$, we get from  \eqref{eqn: cell prob unif conv},
$$ \sup_{j \in \mathbb{N}\cup\{\infty\}} |v_j(R\nu)-\overline{\mu}_j| \leq \delta+ C\omega_\nu(N)^\alpha
 \leq 2\delta+ C\inf_{N>0}\omega_\nu(N)^\alpha. $$
Then fixing such an $R$, since $v_j(R\nu) \to v_\infty(R\nu)$ as $j \to \infty$, for $j>J(\delta,R)$ sufficiently large,
$$|\overline{\mu}_j-\overline{\mu}_\infty| \leq |v_j(R\nu)-\overline{\mu}_j|+|v_j(R\nu)-v_\infty(R\nu)|+|v_\infty(R\nu)-\overline{\mu}_\infty| \leq 5\delta + C\inf_{N>0}\omega_\nu(N)^\alpha.$$
Since $\delta$ was arbitrary we conclude.

\end{proof}

Now let us fix our boundary data $\psi \in C^{\alpha}(\real^n)$ which is $\integer^n$ periodic and $F(M,y)$ our operator satisfying (F1)-(F4).  We will call $\overline{\mu}(\nu) = \overline{\mu}(\psi,F,\nu)$ which is defined in Lemma \ref{lem: ID BC} for all directions $\nu \in S^{n-1} \setminus \real \integer^n$.  We now consider the continuity properties of $\overline{\mu}$ with respect to the normal direction $\nu$.  The goal is to show that $\overline{\mu}$ is continuous at every irrational direction.
\begin{lem}\label{lem: continuity of BC}
Let $\nu \in S^{n-1}\setminus \real\integer^n$.  Then $\overline{\mu}$ satisfies the following continuity estimate for any $N>1$ and $|\nu'-\nu|\leq \eta(\nu,N)$ sufficiently small,
$$ |\overline{\mu}(\nu')-\overline{\mu}(\nu)| \leq C(n,\lambda,\Lambda)(|\psi|_{C^\alpha}+\osc \psi)(N^{\frac{\alpha}{1+\alpha}}|\nu-\nu'|^{\frac{\alpha}{2+\alpha}}+\omega_\nu(N)^\alpha).$$

\end{lem}
\begin{rem}
In particular this estimate implies, since $\omega_\nu(N) \to 0$ as $N \to \infty$ for irrational directions, that $\overline{\mu}(\nu)$ is continuous at every irrational direction.  At every rational direction we have an upper bound on the limiting oscillation of $\overline{\mu}$ in small neighborhoods of $\nu$.  In fact using the methods of \cite{CKL12,ChoiKim12} when the operator $F(M,y)$ homogenizes to a rotationally invariant operator one can show that $\overline{\mu}$ has a unique continuous extension from $S^{n-1} \setminus \real\integer^n$ to all of $S^{n-1}$.  In general though, numerical experiments suggest that there are non rotationally invariant operators for which $\overline{\mu}$ will not extend continuously to the rational directions \cite{CharlieInwon}.  In that case the above continuity estimate may be the best possible.
\end{rem}

\begin{proof}
As before, without loss of generality (due to Lemma \ref{lem: mu properties}) we may assume that $|\psi|_{C^\alpha} \leq 1$ and that $\|\psi\|_\infty = \tfrac{1}{2}\osc \psi $.

\medskip

Fix $\nu \in S^{n-1} \setminus \real\integer^n$ and let $\nu' \in S^{n-1}$ be any direction with $|\nu'-\nu| \leq \eta<1/2$, let $L,R$ with $L>>R$ and $\eta L <<1 $.  The choices of $\eta$, $L$ and $R$ will be made more explicit shortly.  The idea is to explicitly compare the solutions of \eqref{eqn: cell prob} in $P(0,\nu)$ and $P(0,\nu')$ and use the identification of the homogenized boundary condition Lemma \ref{lem: ID BC}.  

\medskip

Let us call $x' = x-(x\cdot \nu)\nu$, then we have,
$$ \{ x \cdot \nu >|\nu'-\nu|L\} \subset P(0,\nu') \cap \{|x'| \leq L \}. $$
Moreover, on $\{x\cdot \nu = |\nu'-\nu|L\}$ letting $a(x') = -\frac{x'\cdot \nu' }{\nu\cdot \nu'}$ so that $x'+a(x')\nu \in \partial P(0,\nu')$ we have that $|a(x')| \leq \frac{1}{2}L|\nu'-\nu|$ since $\eta<1/2$ and by using the $C^{\alpha}$ estimates up to the boundary,
\begin{align*}
|v_{\nu}(x)-v_{\nu'}(x)| &\leq |v_\nu(x)-v_\nu(x')| +|\psi(x')-\psi(x'+a(x')\nu)|+ |v_{\nu'}(x'+a(x')\nu)-v_{\nu'}(x)| \\
& \leq C(n,\lambda,\Lambda)L^\alpha|\nu'-\nu|^\alpha.
\end{align*}
Then we apply the localized comparison estimate of Lemma \ref{lem: localization} in the domain $(P(0,\nu)+ |\nu'-\nu|L\nu) \cap P(0,\nu')$ to get,
$$ |v_{\nu'}-v_{\nu}| \leq C(L^\alpha|\nu'-\nu|^\alpha+\tfrac{R}{L}) \ \hbox{ in } \  \{|x'| \leq R\} \times \{ 0 \leq x \cdot \nu \leq R\}.$$
Then by Lemma \ref{lem: ID BC} for every $1<N< c(n)R$,
\begin{align}
|\overline{\mu}(\nu)-\overline{\mu}(\nu')| &\leq |v_{\nu}(R\nu)-\overline{\mu}(\nu)|+C(L^\alpha|\nu'-\nu|^\alpha+\tfrac{R}{L})+|v_{\nu'}(R\nu)-\overline{\mu}(\nu')| \nonumber \\
& \leq C(L^\alpha|\nu'-\nu|^\alpha+\tfrac{R}{L}+(\tfrac{N}{R})^{\alpha}+\omega_\nu(N)^{\alpha}+\omega_{\nu'}(N)^{\alpha}). \label{eqn: est R L}
\end{align}
Now in order to make the first three terms on the right hand side above of comparable size we choose,
$$ R = N|\nu-\nu'|^{-\frac{1}{2+\alpha}} \ \hbox{ and } \ L = N^{\frac{1}{1+\alpha}}|\nu-\nu'|^{-\frac{1+\alpha(2+\alpha)}{(1+\alpha)(2+\alpha)}}. $$
So for any $N>1$ choose $\eta$ sufficiently small so that $ R\geq C(n) N$, that is,
$$ 0<\eta \leq C(n)^{-(2+\alpha)}. $$
Then based on the continuity properties of $\omega_\cdot(N)$ with fixed $N$ from Section \ref{subsec: Rational and Irrational Directions} choose $\eta$ smaller if necessary so that
$$ \omega_{\nu'}(N) \leq2 \omega_\nu (N) \ \hbox{ when } \ |\nu-\nu'| \leq \eta(\nu,N). $$
Then we get, by plugging in our choice of $R$, $L$ for $|\nu-\nu'| \leq \eta(\nu,N)$ into \eqref{eqn: est R L},
$$ |\overline{\mu}(\nu)-\overline{\mu}(\nu')| \leq C(n,\lambda,\Lambda)(N^{\frac{\alpha}{1+\alpha}}|\nu-\nu'|^{\frac{\alpha}{2+\alpha}}+\omega_\nu(N)^\alpha),$$
which was the desired result.
\end{proof}

\subsection{The Nonlinearity of the Homogenized Boundary Condition} \label{subsec: The Nonlinearity of the Homogenized Boundary Condition}
In the case when the operator $F(M,y) = \Tr (A(y)M)$ is linear the homogenized boundary condition $\overline{\mu}(\psi)$ can be identified explicitly as the average of $\psi$,
\begin{equation}\label{eqn: the average}
\overline{\mu}(\psi,F,\nu) = \int_{[0,1)^n} \psi(y) dy. 
\end{equation}
Note that in particular this value is independent of $\nu \in S^{n-1}\setminus \real \integer^n$ and of the operator $F$.  The proof is a straightforward consequence of the Riesz Representation Theorem in combination with the properties given in Lemma \ref{lem: mu properties}.   We state this in slightly greater generality in the following Lemma,
\begin{lem}\label{lem: linear id}
If $F$ satisfying (F1)-(F4) is concave then for all irrational directions $\nu$ and all $\psi \in C^{\alpha}_{per}(\real^n)$,
$$ \overline{\mu}(\psi,F,\nu)  \geq \int_{[0,1)^n} \psi(y) dy.$$
The opposite inequality holds in case $F$ is convex.  More precisely we argue that if $\overline{\mu}$ is linear in $\psi$ then it is just the average of $\psi$ over the unit periodicity cell.  In particular this is the case when $F$ is a linear operator.
\end{lem}
  
\begin{proof}
We will just show \eqref{eqn: the average}.  Then the result of the Lemma can be derived by noting that $F$ concave (convex) with $F(0,x) = 0$ can be written as a supremum (infimum) of linear homogenous operators and using the monotonicity of $\overline{\mu}$ with respect to $F$.  By Lemma \ref{lem: mu properties} we see that $\overline{\mu} : C^{\alpha}(\mathbb{T}^n) \to \real$ is linear and satisfies,
$$ |\overline{\mu}(\psi)| \leq \|\psi\|_\infty. $$
Therefore $\overline{\mu}$ extends to a bounded linear operator (with norm $1$) on the space $C(\mathbb{T}^n)$, and so, by the Riesz Representation Theorem, there is a Radon measure $\sigma$ on $\mathbb{T}^n$ so that,
\begin{equation*}\label{eqn: linear mu}
 \overline{\mu}(\psi) = \int_{\mathbb{T}^n} \psi(y) d\sigma(y). 
 \end{equation*}
Then, from the translation invariance along with the fact that $\sigma(1) =1$, we derive that $\sigma$ must be the uniform measure.
\end{proof}

Clearly this identification of $\overline{\mu}$, \eqref{eqn: the average}, relies strongly on the linearity of the operator $F$.  Nonetheless one may imagine, since $\overline{\mu}$ is independent of the interior operator in the linear case, that this may carry over to the nonlinear case as well.  This is not the case.  We give a simple example which shows the effect of the nonlinearity:
\begin{ex}\label{ex: nonlinearity}
In dimension $n=2$ let us define, for $\Lambda>1$ and $\nu, \eta$ orthogonal irrational directions, the following fully nonlinear and concave operator,
$$ F(D^2u) = \min\{ -\Delta u, -u_{\eta\eta}-\Lambda u_{\nu\nu}\}.$$
We will take our boundary data to be $\psi(x) = \cos(2\pi x_1)$ and consider the solution $u$ of the following cell problem,
\begin{equation}\label{eqn: example hom}
\left\{\begin{array}{lll}
F(D^2u) = 0 & \hbox{ in } & P(\nu) \\ \\
u(x) = \cos(2\pi x_1) & \hbox{ on } & \partial P(\nu).
\end{array}\right.
\end{equation}
Note that $\int_{[0,1)^n} \psi(x) dx = 0$ would be the homogenized boundary condition associated with $\psi$ if $F$ were linear.  In this case we can construct explicit solutions for both of the operators in the definition of $F$ using separation of variables.  In particular we consider,
$$ v_1(x) = e^{-x \cdot \nu}\cos(2\pi \eta_1(x \cdot \eta)) \ \hbox{ and } \ v_2(x) = e^{-x \cdot \nu/\Lambda^{1/2}}\cos(2\pi \eta_1(x \cdot \eta)), $$
which solve \eqref{eqn: example hom} with the interior operator replaced by $\Delta u$ and $u_{\eta\eta}+\Lambda u_{\nu\nu}$ respectively.  In particular $v_1$ and $v_2$ are both subsolutions of \eqref{eqn: example hom}.  Therefore, by comparison, we see that
\begin{equation}
 u(x) \geq \max\{v_1(x),v_2(x)\} .
 \end{equation}
 Now let us consider the values of $v_1$ and $v_2$ for example on the hyperplane $\{ x \cdot \nu = \Lambda^{1/4}\}$, here we have
 $$ \max\{v_1(x),v_2(x)\} = e^{-\Lambda^{-1/4}}[\cos(2\pi \eta_1(x \cdot \eta))]_++e^{-\Lambda^{1/4}}[\cos(2\pi \eta_1(x \cdot \eta))]_- .$$
 Notice that this is the restriction to the hyperplane $\{ x \cdot \nu = \Lambda^{1/4}\}$ of the $\integer^n$ periodic function on $\real^n$,
 $$ \psi'(x) = e^{-\Lambda^{-1/4}}[\cos(2\pi (x_1-\Lambda^{1/4}\nu_1))]_++e^{-\Lambda^{1/4}}[\cos(2\pi (x_1-\Lambda^{1/4}\nu_1))]_-. $$
Or in other words,
$$ \psi'(x) =\max\{v_1(x),v_2(x)\} \leq u(x) \ \hbox{ on } \ \partial P(\nu,\Lambda^{1/4}\nu).$$
Therefore using Lemma \ref{lem: linear id} identification of $\overline{\mu}$ for linear operators and the monotonicity of $\overline{\mu}$,
$$ \overline{\mu}(\cos(2 \pi \cdot),F,\nu) \geq \overline{\mu}(\psi',-\Delta,\nu) = \int_{[0,1)^n} \psi'(x) dx \geq \tfrac{1}{\pi}\left(e^{-\Lambda^{-1/4}}-e^{-\Lambda^{1/4}}\right)> 0. $$
As a result $\overline{\mu}$ generically will not be just the average of the boundary data when $F$ is nonlinear.  Also note that in particular, since $\overline\mu$ must be the average if it is linear in $\psi$, it cannot be linear in $\psi$.
\end{ex}

\section{Comparison with Partial Boundary Data}\label{subsec: Comparison with Partial Boundary Data}
Recalling that we do not expect the homogenization to hold for \eqref{eqn: epsilon prob} at points on boundary with rational normal direction we are led to consider whether the comparison principle holds when the ordering on the boundary only holds on a `large' subset.  In particular one is led to consider the following situation, suppose that we have $v$ and $u$ bounded and respectively upper and lower semicontinuous in the closure of a bounded domain $\Omega$ and a subset $\Gamma \subset \partial\Omega$ such that,
\begin{equation}\label{eqn: ordered off gamma}
\left\{\begin{array}{lll}
F(D^2v,x)  \leq F(D^2u,x)& \hbox{ in } & \Omega  \\ \\
\limsup_{y \to x} v(y) \leq \liminf_{y \to x} u(y) & \hbox{ on } &  \partial \Omega \setminus \Gamma.
\end{array}\right.
\end{equation}
Then can we deduce that $v \leq u$ in $\Omega$, or more generally is there an estimate of the form,
$$ \sup_{x \in K \subset\subset \Omega} (v-u)_+ \lesssim_K |\Gamma|, $$ 
where $|\cdot|$ is some measure of the size of $\Gamma$.  Such an estimate, a form of the Alexandroff-Bakelman-Pucci estimate on the boundary, appears to be unknown even for linear $F$ at least when $|\cdot| = \mathcal{H}^{n-1}(\cdot)$ the $n-1$ dimensional Hausdorff measure.  We are able to prove such an estimate only for $|\cdot| = \mathcal{H}^{\beta_0}(\cdot)$ with some $\beta_0 <n-1$, possibly quite small depending on the ratio $\Lambda/\lambda$.

\medskip

In order to prove an estimate of the above form we will use the singular solution that was constructed in \cite{ASS12}.  The following Theorem was proven there for even more general operators,
\begin{thm}\textup{(Theorem 1 from \cite{ASS12})}
For any $1>\gamma>0$ let $K_\gamma =  \{x: x_n > (-1+\gamma)|x|\}$, there exists a unique constant $\beta_0(n,\lambda,\Lambda,\gamma)>\max\{\tfrac{\lambda}{\Lambda}(n-1)-1,0\}$ such that the problem,
\begin{equation}\label{eqn: barrier supersoln}
\left\{
\begin{array}{lll}
-\mathcal{P}^+(D^2\Phi) = 0 & \hbox{ in } & K_\gamma \\
\Phi = 0 & \hbox{ on } & \partial K_\gamma\setminus\{0\}
\end{array}\right.
\end{equation}
has a positive solution $\Phi \in C(\overline{K}_\gamma\setminus\{0\})$ which is $-\beta_0$ homogeneous, that is,
$$ \Phi(x) = |x|^{-\beta_0}\Phi(x/|x|). $$
\end{thm}
Note that by the Evans-Krylov Theorem $\Phi$ is $C^{2,\alpha}$ on the interior of $K_\gamma$ for some $\alpha \in (0,1)$, see \cite{CaffarelliCabre95}.  For $\eta \in S^{n-1}$ let $\mathcal{O}_\eta$ be any rotation sending $\eta$ to $e_n$ and then we call, 
\begin{equation}\label{eqn: ASS barrier}
\Phi_\eta(x) = \Phi(\mathcal{O}_\eta x).
\end{equation}

\medskip

In order to use the above barriers we will make the assumption that the domain $\Omega$ has the following exterior cone condition for some $0<\gamma \leq 1$,
\begin{DEF}
Let $0<\gamma < 1$, then a bounded domain $\Omega \subset \real^n$ has the \textit{strict $\gamma$ exterior cone condition} if there exists $\rho>0$ and $\gamma'>\gamma$ such that for every $x \in \partial\Omega$ there is a direction $\eta_x$ with,
$$ \{ y: (y-x)\cdot \eta_x \leq (-1+\gamma')|y-x|\} \cap B_\rho(x)\subset \real^n \setminus \Omega.$$
\end{DEF}

\medskip

\noindent Of course convex domains have the $\gamma$ exterior cone condition for every $\gamma<1$ since they have a supporting hyperplane at every boundary point.  More generally any Lipschitz domain also has the strict $\gamma$ exterior cone condition for some $\gamma$ depending on the Lipschitz constant of $\partial \Omega$.

\medskip

We will also use the notion Hausdorff dimension below, so we will recall the definition.  The $d$-dimensional Hausdorff measure is defined in the following way for a subset $E \subset \real^n$,
$$ \mathcal{H}^d(E) := \sup_{\delta > 0}  \ \inf \left\{ \sum_{j=1}^\infty r_j^d : \exists x_j \in \real^n \hbox{ s.t. } E \subseteq \bigcup_{j\geq1} B(x_j,r_j)  \hbox{ and } r_j \leq \delta\right\}.$$
It is standard to check that $\mathcal{H}^d(E)$ is decreasing in $d$ and $\mathcal{H}^d(E)=0$ for $d>n$, this motivates us to define the Hausdorff dimension as,
$$ \dim_{\mathcal{H}} E := \inf \  \{ d \geq 0 : \mathcal{H}^d(E) = 0 \}. $$

\begin{thm}\label{thm: ABP type}
Let $F$ satisfying (F1) and (F2).  Suppose that $\Omega$ has the strict $\gamma$ exterior cone condition for some $0<\gamma < 1$ and let $\beta_0(n,\lambda,\Lambda,\gamma)$ as above from Theorem 1 of \cite{ASS12}.  Let $\Gamma \subset \partial\Omega$ such that $\mathcal{H}^{\beta_0}(\Gamma) =  0$, in particular this is the case if $\dim_{\mathcal{H}}\Gamma < \beta_0$.  Then if $u$ and $v$ bounded and respectively upper and lower semicontinuous on the closure of $\Omega$ satisfy,
\begin{equation}\label{eqn: ordered off gamma}
\left\{\begin{array}{lll}
F(D^2v,x) \leq F(D^2u,x)& \hbox{ in } & \Omega  \\ \\
\limsup_{y \to x} v(y) \leq \liminf_{y \to x} u(y) & \hbox{ on } &  \partial \Omega \setminus \Gamma
\end{array}\right.
\end{equation}
then $u \leq v$ in $\Omega$.
\end{thm}
\begin{rem}
This Lemma is also true in the gradient dependent case as long as the operator is globally Lipschitz in the gradient argument.  The only alteration is that the best estimate of $\beta_0$.  In fact Theorem 1 of \cite{ASS12} gives the existence of a singular solution in the gradient dependent case as well.
\end{rem}
\begin{proof}
  For simplicity we assume that the $\gamma$ exterior cone condition holds globally, that is $\rho$ from the definition of the $\gamma$ exterior cone condition is equal to $+\infty$.  A localization procedure using Lemma \ref{lem: localization} will cover the general case. Let $w = u-v$ by the Comparison Theorem \ref{thm: comparison pucci} $w$ is a bounded subsolution of
  \begin{equation}\label{eqn: ordered off gamma pucci}
\left\{\begin{array}{lll}
-\mathcal{P}^+(D^2w) \leq 0 & \hbox{ in } & \Omega  \\ \\
\limsup_{y \to x} w(y) \leq 0 & \hbox{ on } &  \partial \Omega \setminus \Gamma.
\end{array}\right.
\end{equation}
   Let $\Phi$ be the singular solution from Theorem 1 of \cite{ASS12} corresponding to $\gamma$, then call $m = \min \{ \Phi(\theta): \theta \in K_{\gamma'} \cap S^{n-1}\}>0$ and take $M> m^{-1}\|w\|_\infty$.  Then from the definition of Hausdorff measure, $\mathcal{H}^{\beta_0}(\Gamma) = 0$ implies that for each $\delta>0$ there exists $B(x_j,r_j)$ for $j \in \mathbb{N}$, an open covering of $\Gamma$ by  balls, so that $ \sum_{j>0} r_j^{\beta_0} < \delta^2/M $. For each $x_j$ let $\eta_j$ be from the $\gamma$ exterior cone condition for $\Omega$ and consider the following barrier, smooth in $\Omega$,
$$ \phi_\delta (x) = \delta+ \sum_{j>0} Mr_j^{\beta_0} \Phi_{\eta_j}(x-x_j). $$
We have the following convergence of the $\phi_\delta \searrow 0$,
$$ 0 \leq \phi_\delta(x) \leq 2\delta \ \hbox{ for } \ d(x,\partial \Omega) \gtrsim \delta^{1/\beta}.$$
Let us show that $\phi_\delta$ is a strict smooth supersolution on the interior of $\Omega$, since $\{y: (y-x_j)\cdot \eta_j \geq(-1+\gamma')|y-x_j|\} \subset \Omega^C$ for each $j$ along with \eqref{eqn: barrier supersoln},
$$ -\mathcal{P}^+(\phi_\delta) \geq  \sum_{j>0} -Mr_j^\beta\mathcal{P}^+(D^2\Phi_{\eta_j}(x-x_j)) \geq 0.$$
Let us show that for each $ x \in \partial \Omega$,
$$ \limsup _{y \to x }  \ (w-\phi_\delta)(y) \leq -\delta.$$
This is true by the assumption for $x \in \partial \Omega \setminus \Gamma$ since $\phi_\delta \geq \delta$ on $\partial \Omega$.  For $x \in \Gamma$ we have that $x \in B(x_j,r_j)$ for some $j>0$, so that, as above, for $ y \in \Omega$ with $|y-x| < d(x, \partial B(x_j,r_j))$,
$$ \phi_\delta(y) \geq \delta+Mr_j^\beta\Phi_{\eta_j}(y-x_j) \geq \delta+Mr_j^\beta mr_j^{-\beta} \geq \delta+\|w\|_\infty.$$
  This proves the ordering of the Dirichlet data.  Then since $\phi_\delta$ is smooth in the interior of $\Omega$ simply by the definition of viscosity solution we have that $w \leq \phi_\delta$ in $\Omega$.  Letting $\delta \to 0$ gives the result.

\end{proof}

\section{Homogenization in General Domains}\label{sec: Homogenization in General Domains}
Now we consider the homogenization of \eqref{eqn: epsilon prob} in general domains $\Omega$.  
We show that \eqref{eqn: epsilon prob} homogenizes when $\partial\Omega$ does not have too many rational boundary points in the sense of Theorem \ref{thm: ABP type}.  Unlike in the situation when $\Omega$ is a half space which was considered in \cite{BarlesMironescu12} it is now necessary to make use of the continuity properties of the homogenized boundary condition $\overline{\mu}$ with respect to $\nu \in S^{n-1}$.  The first reason is in order to guarantee that the homogenized problem,
\begin{equation}
\left\{\begin{array}{lll}
\overline{F}(D^2\overline{u},x) = 0 & \hbox{ in } & \Omega \\ \\
\overline{u}(x) = \overline{g}(x) & \hbox{ on } & \partial \Omega.
\end{array}\right.
\end{equation}
has comparison.  Here $\overline{F}$ is given by Theorem \ref{thm: interior hom} and we will identify $\overline{g}(x)$ by the results of the previous section.  This is where we will use the result of the previous section on comparison with partial boundary data.  

\medskip

In order to blow up at boundary points and use the results for the cell problem in a locally uniform way we will again need the stability of the rate of averaging at irrational directions. We make this precise in Lemma \ref{lem: bdry layer general} which is based on Lemma 3.1 in \cite{BarlesMironescu12}, but first we make the following note on the behavior of $F_\e$ under the blow up rescaling:
\begin{lem}\label{lem: scale one op}
The following convergence holds:
\begin{equation*}
\varepsilon^2F_\varepsilon( \varepsilon^{-2}M,x,y) \to F^x(M,y) \ \hbox{ as } \ \varepsilon \to 0,
\end{equation*}
locally uniformly in $M$ and uniformly in $(x,y)$, with $F^x(M,y)$ satisfying assumptions (F1)-(F4).  In particular we also have that if $x_j \to x$ then $F^x_j \to F^x$ locally uniformly in $M$ and uniformly in $y$.  
\end{lem}

\begin{proof}
By our assumptions on $F_\e$ in Section \ref{subsec: Assumptions on the Operators}, there exists a collection of matrices, $\lambda I \leq A^{\alpha\beta}(x,y) \leq \Lambda I$, $A^{\alpha,\beta}$ have uniformly bounded Lipschitz norm and are periodic in their $y$ variable, also there exists a function $f \in C^{0,1}(\real^n \times \real^n)$ periodic in it's second variable so that,
$$F_\varepsilon(M,x,y) = \inf_{\beta \in \mathcal{B}} \sup_{\alpha \in \mathcal{A}}\left[-\Tr(A^{\alpha\beta}(x,y)M)\right]+f(x,y).$$
This form of the operators makes the claim clear, and in this situation,
\begin{equation}
 F^x(M,y) = \inf_{\beta \in \mathcal{B}} \sup_{\alpha \in \mathcal{A}}\left[-\Tr(A^{\alpha\beta}(x,y)M)\right].
 \end{equation}
\end{proof}

\begin{lem}\label{lem: bdry layer general}
 \textup{(Behavior of $u^\varepsilon$ away from the boundary layer)}
Suppose that $\Omega$ is a bounded domain with $\partial \Omega$ being $C^2$.  Let $x_0 \in \partial \Omega$ such that the interior normal $ \nu_{x_0}$ is an irrational direction.   For any $\delta>0$ there exists $R_0$ and $\eta>0$ such that for $R \geq R_0$ and for any sequence $x_\varepsilon \to x$ with $|x-x_0| \leq \eta$,
$$ \limsup_{\varepsilon \to 0} |u^\varepsilon(x_\varepsilon+\varepsilon R\nu_{x_0})-\overline{g}(x)| \leq \delta,$$
(in particular the quantity in the $\limsup$ is defined for $\e$ sufficiently small) where the homogenized boundary data $\overline{g}$ is given by,
$$ \overline{g}(x) := \overline{\mu}(g(x,\cdot),F^x,\nu_{x}). $$
\end{lem}

\begin{proof}
First we explain how to choose $\eta$.  Given $\delta>0$, and $\nu_{x_0}$ an irrational direction, there exists $N_0>1$ sufficiently large so that $
\omega_{\nu_{x_0}}(N_0) \lesssim \delta^{1/\alpha}$.  Recalling that $\omega_\cdot$ is continuous at irrational directions and that $\nu_x$ is continuous in $x$ since $\Omega$ is assumed to be a $C^2$ domain, let $\eta$ be sufficiently small so that 
$$ |x-x_0| \leq \eta \ \hbox{ implies } \ \omega_{\nu_x} (N) \lesssim \delta^{1/\alpha}.$$
Now fix $x \in \partial \Omega$ be fixed as above with $|x-x_0| \leq \eta$ as above and for every $\varepsilon>0$ define
$$ v^\varepsilon(y) := u^\varepsilon(x+\varepsilon y). $$
Let us call $\Omega^\e = \e^{-1}(\Omega-x)$.  The $v^\varepsilon$ satisfy the following equation in the viscosity sense,
\begin{equation}\label{eqn: scale one prob}
\left\{
\begin{array}{lll}
F_\varepsilon(\varepsilon^{-2}D^2v^\varepsilon,x+\varepsilon y,\varepsilon^{-1}x+y) = 0 & \hbox{ for } & y \in  \Omega^\e\\ \\
v^\varepsilon(y) = g(x+\varepsilon y,\varepsilon^{-1}x+y) & \hbox{ for } & y \in \partial \Omega^\e.
\end{array}\right.
\end{equation}
The $\Omega^\e \to P_x=P(0,\nu_x)$ in Hausdorff distance when restricted to compact sets.  Now we use the compactness afforded by the periodic setting.  Let us call
$$ [0,1)^n \ni \tau_\e = \e^{-1} x \mod \integer^n,$$
 i.e. $\varepsilon^{-1}x =\tau_\varepsilon+z_\varepsilon$ with $z_\varepsilon \in \integer^n$.  We may replace $\varepsilon^{-1}x$ in \eqref{eqn: scale one prob} above by $\tau_\varepsilon$ using the $\integer^n$ periodicity of $F_\varepsilon$ and $g$.  Then taking a subsequence of the $\varepsilon$ the $\tau_{\varepsilon} \to \tau \in [0,1]^n$.  Noting that $ -\|g\|_{L^\infty(\real^n\times \real^n)} \leq v^\varepsilon \leq \|g\|_{L^\infty(\real^n\times \real^n)} $ we may take the upper and lower relaxed limits of $v^{\varepsilon}$ along the above subsequence,
\begin{equation*}
\begin{array}{lll}
v^{*,\tau}(y) = \left.\limsup\right.^* v^\varepsilon(z) & \hbox{ and } & v_{*,\tau}(y) = \left.\liminf\right._*  v_\varepsilon(z) .
\end{array}
\end{equation*}
The upper and lower half-relaxed limits are defined for each $y \in \overline{P_x}$.  Let us show that these are respectively sub and supersolutions of the limiting scale one equation,
\begin{equation}\label{eqn: scale one prob}
\left\{
\begin{array}{lll}
F^x(D^2v_\tau,y+\tau) = 0 & \hbox{ for } & y \in P_x \\ \\
v_\tau(y) = g(x,\tau+y) & \hbox{ for } & y \in \partial P_x.
\end{array}\right.
\end{equation}
The interior sub/supersolution property follows from the stability of viscosity solutions under uniform convergence, Lemma \ref{lem: uniform stability}.  The condition on the boundary can be shown in the following way.  Fix $y \in \partial P_x$ and $\delta>0$.  For $\e$ sufficiently small, from the local Hausdorff distance convergence of $\partial\Omega^\e \to \partial P_x$, there exists $y' \in \partial \Omega_\e \cap B(y,\delta)$.  Now for any $z \in \Omega^\e \cap B(y,\delta)$, by the estimates up to the boundary Lemma \ref{lem: bdry estimates},
\begin{align*} 
|v^\e(z)-g(x,\tau+y)| &\leq C\|g(x+\e \cdot, \tau_\e+\cdot)\|_{C^\alpha}|z-y'|^\alpha+|g(x+\e y', \tau_\e+y')- g(x, \tau+y)| \\
& \leq  C(\|g\|_{C^\alpha})(\delta^\alpha+\e^\alpha|y'|^\alpha+|\tau_\e-\tau|^\alpha).
\end{align*}
Then taking the supremum over such $z$ and letting $\e \to 0$ along any subsequence where $\tau_\e \to \tau$,
$$|v^{*,\tau}(y)-g(x,\tau+y)| \leq \lim_{\e \to 0} \sup_{z \in \Omega^\e \cap B(y,\delta)} |v^\e(z)-g(x,\tau+y)| \leq C\|g\|_{C^\alpha} \delta^\alpha.$$
Letting $\delta \to 0$ proves that $v^{*,\tau}$ achieves the boundary data in \eqref{eqn: scale one prob}.  The same argument works for $v_{*,\tau}$.

\medskip

Since the equation \eqref{eqn: scale one prob} has comparison we derive $v^{*,\tau} \leq v_{*,\tau}$, of course the opposite inequality is always true so we get the local uniform convergence in $P_x$ along subsequences of $v^\varepsilon \to v_\tau$ for some $\tau \in [0,1]^n$ depending on the subsequence.  The point is that the limiting behavior of $v_\tau(y+R\nu)$ will not depend on $\tau \in [0,1]^n$.

\medskip

Let us consider $\widetilde{v}_\tau = v_\tau (y-\tau)$, these satisfy,
 \begin{equation}
\left\{
\begin{array}{lll}
F^x(D^2\widetilde{v}_\tau,y) = 0 & \hbox{ for } & y \in P_x+(\tau \cdot \nu_x)\nu_x \\ \\
\widetilde{v}_\tau(y) = g(x,y) & \hbox{ for } & y \in \partial P_x+(\tau \cdot \nu_x)\nu_x.
\end{array}\right.
\end{equation}
Then by Lemma \ref{lem: ID BC}, for any $N>1$,
$$\sup_{ \tau \in [0,1]^n}\sup_{y \in  P_x +(\tau \cdot \nu_x)\nu_x} |\widetilde{v}_\tau(y+R\nu_x)-\overline{\mu}(g(x,\cdot),F^x,\nu_x)| \lesssim \omega_{\nu_x}(N)^\alpha+(N/R)^{\alpha}.$$
Given $\delta>0$ fixing $N=N_0$ as in our choice of $\eta$ above and letting $R_0>0$ sufficiently large we get,
$$\sup_{ \tau \in [0,1]^n}\sup_{y \in \partial P_x +(\tau \cdot \nu_x)\nu_x} |\widetilde{v}_\tau(y+R\nu_x)-\overline{g}(x)| \leq \delta. $$
Then for $R\geq 2R_0$, $\eta>0$ sufficiently small by the continuity of $\nu_x$ so that $|\nu_x-\nu_{x_0}| \leq 1/2$, and each $\tau \in [0,1]^n$,
$$|u^\varepsilon(x+\varepsilon R\nu_{x_0})-\overline{g}(x)| \leq |u^\varepsilon(x+\varepsilon R\nu_{x_0})-v_\tau(R\nu_{x_0})|+|v_\tau(R\nu_{x_0})-\overline{g}(x)| \leq |u^\varepsilon(x+\varepsilon R\nu_{x_0})-v_\tau(R\nu)|+\delta,$$
and infimizing over $\tau$ and then taking $\limsup$ on both sides,
$$\limsup_{\varepsilon \to 0}|u^\varepsilon(x+\varepsilon R\nu_{x_0})-\overline{g}(x)| \leq \delta+\limsup_{\varepsilon \to 0}\inf_{\tau \in [0,1]^n}|u^\varepsilon(x+\varepsilon R\nu_{x_0})-v_\tau(R\nu_{x_0})| = \delta.$$
Since the convergence of the $u^\varepsilon$ is locally uniform along subsequences we may replace $x$ by any sequence $x_\varepsilon \to x$ above.
\end{proof}

  From Lemma \ref{lem: bdry layer general} we will be able to deduce that homogenization of the boundary condition holds near boundary points with irrational directions.  In general one cannot blow up to get an estimate of the form in Lemma \ref{lem: bdry layer general} near boundary points with rational normal directions.  However we have the comparison principle with partial boundary data from Theorem \ref{thm: ABP type}.  This Theorem gives a condition under which satisfying the Dirichlet data only at boundary points with irrational normal directions is enough to deduce uniqueness.

\medskip

Now let $\Omega \subset \real^n$ be a bounded domain with $C^2$ boundary which has a strict $\gamma$ exterior cone condition for some $0<\gamma < 1$.  Recalling that $\nu_\cdot : \partial \Omega \to S^{n-1}$ is the inward unit normal vector field to $\partial \Omega$ we call the set of rational boundary points,
$$ \Gamma = \Gamma(\Omega) := \{ x\in \partial \Omega : \nu_x \in \real \integer^n\} .$$
Let us suppose that $\Omega$ has no flat boundary pieces in any rational direction of too large Hausdorff dimension, precisely we assume,
$$ \mathcal{H}^{\beta_0}(\Gamma) = 0, $$
where $n-1>\beta_0(n,\lambda,\Lambda,\gamma)>\max\{0,\tfrac{\lambda}{\Lambda}(n-1)-1\}$ is the singular solution exponent from Theorem 1 in \cite{ASS12}.  In particular this is always true of uniformly convex domains which have an exterior half space at every boundary point and for which $\Gamma$ will be countable.  More generally any bounded $C^2$ domain such that the set of rational boundary points where $D_\tau \nu_x$ (the tangential gradient) is degenerate has Hausdorff dimension smaller than $\beta_0$ will satisfy this condition.

\begin{thm}\label{thm: dirichlet hom}
Suppose that $\Omega$ satisfies the hypothesis given above.  Let $\overline{g}(x) = \overline{\mu}(g(x,\cdot),F^x,\nu_x)$, then $\overline{g}: \partial \Omega \to \real$ is continuous at every point of $\partial \Omega \setminus \Gamma$, and the solutions $u^\varepsilon$ of \eqref{eqn: epsilon prob} converge locally uniformly in $\Omega$ to the unique, in the sense of Theorem \ref{thm: ABP type}, solution given by Perron's method $\overline{u}$ of
\begin{equation}\label{eqn: homogenized eqn gen}
\left\{
\begin{array}{lll}
\overline{F}(D^2\overline{u},x) = 0 & \hbox{ in } & \Omega \\
\overline{u}(x) = \overline{g}(x) & \hbox{ on } & \partial \Omega \setminus \Gamma.
\end{array}\right.
\end{equation}
\end{thm}

\begin{proof}
Let us first note the continuity of $\overline{g}$ at points of $\partial\Omega \setminus\Gamma$.  Let $x_j \in \partial \Omega $ converging to $x \in \partial \Omega \setminus \Gamma$, then $\|g(x_j,\cdot)-g(x,\cdot)\|_{L^\infty} \leq C|x_n-x|^\alpha$, $\nu_{x_n} \to \nu_x$ since $\partial \Omega$ is $C^2$ and $F^{x_j}(M,y)$ converges to $F^x(M,y)$ locally uniformly in $M$ and uniformly in $y$ as $j \to \infty$.  Then from the continuity of $\overline{\mu}$ noted in lemma \ref{lem: mu properties} (i) and (iii) and in lemma \ref{lem: continuity of BC} we have the convergence, 
\begin{equation}\label{eqn: cont final}
\overline{g}(x_j)=\overline{\mu}(g(x_j,\cdot),F^{x_j},\nu_{x_j}) \to \overline{\mu}(g(x,\cdot),F^x,\nu_x)=\overline{g}(x) \ \hbox{ as } \ j \to \infty.
\end{equation}
This shows the claimed continuity of $\overline{g}$.

\medskip

The main tool to show the homogenization is Lemma \ref{lem: bdry layer general}.  Fixing $x_0 \in \partial\Omega \setminus S$, $\nu_{x_0}$ the interior normal to $\partial \Omega$ at $x_0$, $\delta>0$ and $R\geq R_0(\delta)$ and $\eta = \eta(\delta)>0$ from Lemma \ref{lem: bdry layer general} we consider,
$$\widetilde{u}^\varepsilon(x) := u^\varepsilon(x+\varepsilon R\nu_0) \ \hbox{ defined for } \ x \in \Omega - \e R \nu_0.$$
Since the $\widetilde{u}^\varepsilon$ are uniformly bounded (by comparison with $\pm \|g\|_\infty$) the upper and lower relaxed limits are defined and finite in $\overline{\Omega}$,
\begin{equation*}
\begin{array}{lll}
\widetilde{u}^* := \left.\limsup\right.^* \widetilde{u}^\varepsilon & \hbox{ and } & \widetilde{u}_* := \left.\liminf\right._* \widetilde{u}^\varepsilon .
\end{array}
\end{equation*}
Note that for $x \in \Omega$ the upper (and lower) relaxed limits of $\widetilde{u}^\varepsilon$ and $u^\varepsilon$ (which we call $u^*$ and $u_*$ respectively) agree.  The only dependence on $R$, the manifestation of the boundary layer for $u^\varepsilon$, is on $\partial \Omega$.  By our choice of $R$ and $\eta$ from Lemma \ref{lem: bdry layer general} and the continuity shown above in \eqref{eqn: cont final} there exists $0<\eta'\leq\eta$ sufficiently small so that for $x \in \partial \Omega$ with $|x-x_0| \leq \eta'$,
\begin{equation*}
\begin{array}{lll}
\widetilde{u}^*(x) \leq \overline{g}(x_0)+\delta & \hbox{ and } & \widetilde{u}_*(x) \geq \overline{g}(x_0)-\delta,
\end{array}
\end{equation*} 
and by the interior homogenization result for $x \in \Omega$ we get the inequalities,
\begin{equation*}
\begin{array}{lll}
\overline{F}(D^2\widetilde{u}^*,x) \leq 0 & \hbox{ and } & \overline{F}(D^2\widetilde{u}_*,x) \geq 0.
\end{array}
\end{equation*}
 Let $\phi_{\pm}$ be respectively super and subsolution barriers for the domain $\Omega$ at $x_0$.  Then, as in the standard barrier argument, we consider
 $$ \widetilde{\phi}_+=\frac{\|g\|_\infty}{\inf_{|x-x_0| \geq \eta} \phi_+}\phi_++(g(x_0)+ \delta),$$
 which satisfies $\widetilde{\phi}_+ \geq \widetilde{u}^*$ on $\partial\Omega$.  By comparison $\widetilde{\phi}_+ \geq \widetilde{u}^*$ in $\Omega$.  In fact since $\phi_\pm$ can be chosen smooth this is even just by the definition of subsolution.  In any case, since $u^* = \widetilde{u}^*$ in $\Omega$ we get,
 $$ \limsup_{z \to x_0} \ u^*(z)\leq \limsup_{z \to x_0} \ \widetilde{u}^*(z) \leq g(x_0)+\delta. $$
 Making the same argument with $\phi_-$ and $\widetilde{u}_*$ we get,
 $$ \liminf_{z \to x_0} \ u_*(z) \geq g(x_0)-\delta. $$
 Of course now the left hand side above is independent of $\delta$, so we may take $\delta \to 0$ to get,
 $$ g(x_0) \leq \liminf_{z \to x_0} \ u_*(z) \leq \limsup_{z \to x_0} \ u^*(z) \leq g(x_0). $$
 This argument works for every $x_0 \in \partial \Omega \setminus \Gamma$.  Applying the comparison principle with partial boundary data Theorem \ref{thm: ABP type} for $\overline{F}$ to $u^*$ and $u_*$ we get for $x \in \Omega$,
$$  u^*(x)=u_*(x) ,$$
so that the $u^\varepsilon$ converge locally uniformly in $\Omega$ to $\overline{u} = u^*=u_*$.

\end{proof}

\bibliographystyle{plain}
\bibliography{dirichletArticles}
\end{document}